\documentclass[12pt,a4paper]{article}
\usepackage{bbm}
\usepackage{mathrsfs}
\usepackage{graphicx}
\usepackage{amsmath}
\usepackage{amssymb}
\usepackage{amsfonts}
\usepackage{enumerate}
\usepackage{theorem}

\makeatletter
\def\tank#1{\protected@xdef\@thanks{\@thanks
        \protect\footnotetext[0]{#1}}}
\def\bigfoot{

    \@footnotetext}
\makeatother

\newcommand{\ea}{\end{array}}
\newtheorem{theorem}{Theorem}[section]
\newtheorem{proposition}{Proposition}[section]

\newtheorem{lemma}{Lemma}[section]
\newtheorem{definition}{Definition}[section]
\newtheorem{Rem}{Remark}[section]

{\theorembodyfont{\rmfamily}
}

\newenvironment{proof}{Proof.}

\title
{\bf Finite dimensionality of global attractor for the solutions to 3D viscous primitive equations  of large-scale moist atmosphere \thanks{This work was partially
supported by NNSF of China(Grant No. 11401057),   Natural Science Foundation Project of CQ  (Grant No. cstc2016jcyjA0326),
Fundamental Research Funds for the Central Universities(Grant No. 106112015CDJXY100005) and China Scholarship Council (Grant No.201506055003).} }
\author{
 Guoli Zhou
\thanks{ Chongqing University, P.R. China }
\tank{E-mail:zhouguoli736@126.com.}
}
\begin{document}
\maketitle

\begin{abstract}
Under general boundary conditions we  consider the finiteness of  the Hausdorff and fractal dimensions of the global attractor for the strong solution of the 3D moist primitive equations with viscosity.  Firstly, we  obtain time-uniform estimates of the first-order time derivative of the strong solutions in $L^2(\mho)$.  Then,  to prove the finiteness of  the Hausdorff and fractal dimensions of the global attractor,  the  common method is to  obtain the uniform boundedness of the strong solution in $H^2(\mho)$ to establish the squeezing property of the solution operator. But it is difficult to achieve due to the boundary conditions and complicated structure of the 3D moist primitive equations. To overcome the difficulties, we try to use the uniform boundedness of the derivative of the strong solutions with respect to time $t$ in $L^2(\mho)$ to prove the uniform continuity of the global attractor.  Finally, using the uniform continuity of the global attractor we  establish the squeezing property of the solution operator which implies the  finiteness of the Hausdorff and fractal dimensions of the global attractor.
\end{abstract}


\section{Introduction}
The paper is concerned with the 3-dimensional viscous primitive equations in the pressure coordinate system (see\ e.g. $\cite{GH1, LC, LTW1, LTW2}$ and the references therein).
 \begin{eqnarray}
&&\partial_{t}v+\nabla_{v}v+w\partial_{\xi}v+\frac{f}{R_{0}}v^{\bot} +\mathrm{grad} \Phi+ L_{1}v=0,\\
&&\partial_{\xi}\Phi+\frac{bP}{p}(1+aq)T=0, \\
&&\mathrm{div} v+ \partial_{\xi} w=0, \\
&&\partial_{t} T+\nabla_{v}T +w\partial_{\xi}T-\frac{bP}{p}(1+aq)w+L_{2}T=Q_{1} ,\\
&&\partial_{t}q+ \nabla_{v}q +w\partial_{\xi}q+L_{3}q=Q_{2}  .
\end{eqnarray}
The unknowns for the primitive equations are the fluid velocity field $(v,w )=(v_{\theta},v_{\varphi},w )\in \mathbb{R}^{3}$ with  $v=(v_{\theta},v_{\varphi})$ and $ v^{\perp}=(-v_{\varphi}, v_{\theta} ) $ being horizontal, the temperature $T$, $q$ the mixing ratio of water vapor in the air  and the geopotential $\Phi.$
$f=2 \mathrm{cos\theta}$ is the given Coriolis parameter, $Q_{1}$ corresponds to the sum of the heating of the sun and the heat added or removed by condensation or evaporation, $Q_{2}$ represents the amount of water added or removed by condensation or evaporation, $a$ and $b$ are positive constants with $a\approx 0.618$, $R_{0}$ is the Rossby number, $P$ stands for an approximate value of pressure at the surface of the earth, $p_{0}$ is the pressure of the upper atmosphere with $p_{0}>0$ and the variable $\xi$ satisfies $p=(P-p_{0})\xi+p_{0}$ where $0<p_{0}\leq p\leq P.$ The viscosity, the heat and the water vapor diffusion operators $L_{1},\ L_{2}$ and $L_{3}$ are given respectively as the following:
$$L_{i}=-\nu_{i}\Delta-\mu_{i}\partial_{zz}, i=1,2,3. $$
Here the positive constants $\nu_{1}, \mu_{1}$ are the horizontal and vertical viscosity coefficients;  the positive constant $\nu_{2}, \mu_{2}$ are the horizontal and vertical heat diffusivity coefficients; while the positive constant $\nu_{3}, \mu_{3}$ are the horizontal and vertical water vapor diffusivity coefficients. The definitions of $\nabla_{v}v, \Delta v, \Delta T, \Delta q, \nabla_{v}q, \nabla_{v}T, \mathrm{div} v, \mathrm{grad} \Phi $ will be given in section 2.
\par
The space domain of equations: $(1.1)-(1.5)$ is
$$\mho=S^{2}\times (0,1), $$
 where $S^{2}$ is two-dimensional unit sphere. The boundary value conditions are given by
\begin{eqnarray}
\xi=1 (p=P): \partial_{\xi}v=0,\ \ w=0,\ \ \partial_{\xi}T=\alpha_{s}(T_{s}-T),\ \ \partial_{\xi}q=\beta_{s}(q_{s}-q),
\end{eqnarray}
\begin{eqnarray}
\xi=0 (p=p_{0}): \partial_{\xi}v=0,\ \ w=0,\ \ \partial_{\xi}T=0,\ \ \partial_{\xi}q=0,
\end{eqnarray}
where $\alpha_{s}, \beta_{s}$ are positive constants, $T_{s}$ is the given temperature on the surface of the earth, $q_{s}$ is the given mixing ratio of water vapor on the surface of the earth. To simplify the notations, we set $T_{s}=0$ and $q_{s}=0$ without losing any generality. For the case $T_{s}\neq 0$ and $q_{s}\neq 0,$ we can homogenize the boundary value conditions for $T,q;$ see $\cite{GH1}$ for detailed discussion on this issue. Moreover, using $(1.2), (1.3)$ and the boundary conditions $(1.6)-(1.7)$, we have
\begin{eqnarray}
w(t;\theta,\varphi,\xi)=\int_{\xi}^{1}\mathrm{div}\ v( t;\theta,\varphi,\xi')d\xi',
\end{eqnarray}
\begin{eqnarray}
\int_{0}^{1}\mathrm{div}\ v d\xi=0,
\end{eqnarray}
\begin{eqnarray}
\Phi( t;\theta,\varphi,\xi)=\Phi_{s}( t;\theta,\varphi)+\int_{\xi}^{1}\frac{bP}{p}(1+aq)Td\xi',
\end{eqnarray}
where $\Phi_{s}( t;\theta,\varphi)$ is a certain unknown function at the isobaric surface $\xi=1.$ In this article, we assume that the constants $v_{i}=\mu_{i}=1, i=1,2,3. $ For the general case, the results will still be valid.
Then using $(1.8)-(1.10)$, we obtain the following equivalent formulation for system $(1.1)-(1.7)$ with initial condition
\begin{eqnarray}
\partial_{t}v&+&\nabla_{v}v+ \Big{(} \int_{\xi}^{1}\mathrm{div}\ v( t;\theta,\varphi,\xi')d\xi'\Big{)}  \partial_{\xi}v+\frac{f}{R_{0}}v^{\bot} +\mathrm{grad} \Phi_{s}\nonumber\\
&&+\int_{\xi}^{1}\frac{bP}{p}\mathrm{grad}[(1+aq)T]d\xi'  -\Delta v- \partial_{\xi\xi}v =0,
\end{eqnarray}
\begin{eqnarray}
\partial_{t} T&+&\nabla_{v}T +\Big{(} \int_{\xi}^{1}\mathrm{div}\ v( t;\theta,\varphi,\xi')d\xi'\Big{)}\partial_{\xi}T\nonumber\\
&&-\frac{bP}{p}(1+aq)\Big{(} \int_{\xi}^{1}\mathrm{div}\ v( t;\theta,\varphi,\xi')d\xi'\Big{)}-\Delta T- \partial_{\xi\xi}T =Q_{1} ,
\end{eqnarray}
\begin{eqnarray}
\partial_{t} q+\nabla_{v}q +\Big{(} \int_{\xi}^{1}\mathrm{div}\ v( t;\theta,\varphi,\xi')d\xi'\Big{)}\partial_{\xi}q -\Delta q- \partial_{\xi\xi}q =Q_{2} ,
\end{eqnarray}
\begin{eqnarray}
\int_{0}^{1}\mathrm{div}\ vd\xi=0,
\end{eqnarray}
\begin{eqnarray}
\xi=1 : \partial_{\xi}v=0,\ \ w=0,\ \ \partial_{\xi}T=-\alpha_{s}T,\ \ \partial_{\xi}q=-\beta_{s}q,
\end{eqnarray}
\begin{eqnarray}
\xi=0 : \partial_{\xi}v=0,\ \ w=0,\ \ \partial_{\xi}T=0,\ \ \partial_{\xi}q=0,
\end{eqnarray}
\begin{eqnarray}
v(0;\theta,\phi, \xi)=v_{0}(\theta,\phi, \xi ), T(0;\theta,\phi, \xi)=T_{0}(\theta,\phi, \xi ), q(0;\theta,\phi, \xi)=q_{0}(\theta,\phi, \xi ).
\end{eqnarray}
In order to understand the mechanism of long-term weather prediction, one can take advantage of the historical records and numerical computations to detect the future weather. Alternatively, one should also study the long time behavior mathematically for the equations and models governing the motion. The primitive equations represent the classic model for the study of climate and weather prediction, describing the motion of the atmosphere when the hydrostatic assumption is enforced $\cite{ G, Ha, HW, MT, Ri}.$ But the resulting flow or the atmosphere is rich in its organization and complexity (see $\cite{G, Ha,HW}$), the full governing equations are too complicated to be treatable both from the theoretical and the computational side. To overcome this difficulty, some simple numerical models were introduced. The 2-D and 3-D quasi-geostrophic models have been the subject of analytical mathematical study (see e.g., \cite{BB, C, CMT1, CMT2, CW, EM, M, W1, W2, W3} and references therein).
To the best of our knowledge, the mathematical framework of primitive equations was formulated in $\cite{LTW1, LTW2, LTW3},$
where the definitions of weak and strong solutions were given and the existence of weak solution was proven, leaving the uniqueness of weak solution as an open problem for now. Local well-posedness of strong solutions was obtained in $\cite{GMR, TZ}.$ If the domains was thin, the global well-posedness of 3D primitive equations was shown in $\cite{HTZ}.$  Taking advantage of the fact that the pressure is essentially two-dimensional in the primitive equations, global well-posedness of the full three-dimensional case was established in $\cite{CT1}$ and independently in $\cite{Kob1, Kob2}.$ In the subsequent work $\cite{KZ}$ a different proof was developed which allows one to treat non-rectangular domains. Recently, the results were improved in $\cite{CLT1, CLT2, CLT3, CT2}$  by considering the system with partial dissipation, i.e. , with only partial viscosities or only partial diffusion. For the inviscid primitive equations, finite-time blowup was established in $\cite{CINT}.$ To study the long term behavior of primitive equations, the existence of global attractor was established in $\cite{J}$ and dimensions were proven to be finite in $\cite{JT}.$ When moisture is included, an equation for the conservation of water must be added, which is the case in e.g. $\cite{GH1,GH2, LTW1, PTZ}.$ In $\cite{ZHKTZ},$ global well-posedness of quasi-strong and strong solutions was obtained for the primitive equations of atmosphere in presence of vapour saturation.
\par
The understanding of asymptotic behavior of dynamical system is one of the most important topics of modern mathematical physics. One way to solve the problem for dissipative deterministic dynamical system is to consider its global attractor (see its definition in section 2). Thus, in order to capture the dynamical features of moist primitive equations, Guo and Huang in $\cite{GH2}$ proved the existence of universal attractor.  Recently, the authors in $\cite{ZG}$ proved the existence of the global attractor of the strong solutions to the 3D moist primitive equations.
\par
In this paper, we investigate the finiteness of the
the Hausdorff and fractal dimensions of the global attractor for the strong solutions to the 3D moist primitive equations with viscosity. The common method is to obtain uniformly estimates for the strong solution in $H^{2}(\mho)$ space and then use the squeezing property of the solution operator to show our main result. But it is difficult to get the time-uniform estimates for the strong solution in the space $H^{2}(\mho)$. Because, firstly, the boundary conditions in this work is general and consists with the boundary conditions assuring the global well-posedness of the moist primitive equations without any extra conditions which help us to use integration by parts formula  to obtain $a\ priori$ estimates in space with higher regularity.
 Secondly, the structure of the moist primitive equations is even more complicated than the oceanic primitive equations studied in $\cite{CT1}$. For example, in the horizontal momentum equation, there is a challenging term of the gradient of temperature times the mixing ratio of water vapor. The temperature equations and the mixing ratio of water vapor equations also have the similar challenging terms which  present essential difficulties for $a$ $priori$ estimates in the $H^{2}(\mho)$ norm.  To overcome the difficulties, inspired by $\cite{J1}$ we try to use the uniform continuity of the global attractor to prove the squeezing property of the solution operator which implies the finiteness of the
the Hausdorff and fractal dimensions of the global attractor for the strong solutions to the 3D moist primitive equations. This method can be applied to other dissipative equations with physical boundary conditions.
\par
The remaining of the paper is organized as follows. In section $2$, we present the notations and recall some important facts which are crucial to later analysis. Absorbing ball of the first-order time derivatives of the solution in $L^{2}(\mho)$ is obtained in section $3$. Section $4$ is for the uniform continuity and the dimensions of the global attractor. As usual, the positive constants $c$ may change from one line to the next, unless, we give a special declaration.

\section{Preliminaries }
In this section we collect some preliminary results that will be used in the rest of this paper, and we start with the following notations which will be used throughout this work.
Denote
\begin{eqnarray*}
\bar{v}=\int_{0}^{1}vd\xi ,\ \   \tilde{v}=v-\bar{v} .
\end{eqnarray*}
Then we have
\begin{eqnarray}
\nabla \cdot  \bar{v}=0,\ \   \bar{\tilde{v}}=0.
\end{eqnarray}
\par
Now we give the definitions of some differential operators. Firstly,  the natural generalization of the directional derivative on the Euclidean space to the covariant derivative on $S^{2}$ is given as follows.  Let $T,q, \in C^{\infty}(\mho),   \Phi_{s} \in C^{\infty}(S^{2} )$ and
\begin{eqnarray*}
v=v_{\theta}e_{\theta}+v_{\varphi}e_{\varphi},\ \ \ \ \ u=u_{\theta}e_{\theta} +u_{\varphi}e_{\varphi}\ \ \ \ \in C^{\infty}(T\mho| TS^{2}),
\end{eqnarray*}
where $C^{\infty}(T\mho| TS^{2}) $ is the first two components of smooth vector fields on $\mho.$  We define  the covariant derivative of $u,T $ and $q$ with respect to $v$  as follows
\begin{eqnarray*}
\nabla_{v}u=(v_{\theta}\partial_{\theta}u_{\theta}+\frac{v_{\varphi}}{\mathrm{sin}\theta}\partial_{\varphi}u_{\theta}-v_{\varphi}u_{\varphi}\mathrm{cot}\theta )e_{\theta}+(v_{\theta}\partial_{\theta}u_{\varphi}+\frac{v_{\varphi}}{\mathrm{sin}\theta}\partial_{\varphi}u_{\varphi}+v_{\varphi}u_{\theta}\mathrm{cot}\theta )e_{\varphi},
\end{eqnarray*}
\begin{eqnarray*}
\nabla_{v}T= v_{\theta}\partial_{\theta}T+ \frac{v_{\varphi}}{\mathrm{sin}\theta}\partial_{\varphi}T,
\end{eqnarray*}
\begin{eqnarray*}
\nabla_{v}q= v_{\theta}\partial_{\theta}q+ \frac{v_{\varphi}}{\mathrm{sin}\theta}\partial_{\varphi}q.
\end{eqnarray*}
We give the definition of the horizontal gradient $\nabla = \mathrm{grad}$ for $T$ and $\Phi_{s}$ on $S^{2}$ by
\begin{eqnarray*}
\nabla T= \mathrm{grad} T= (\partial_{\theta}T) e_{\theta} +\frac{1}{\mathrm{sin}\theta}(\partial_{\varphi}T)e_{\varphi},
\end{eqnarray*}
\begin{eqnarray*}
\nabla \Phi_{s}= \mathrm{grad} \Phi_{s}=( \partial_{\theta}\Phi_{s} )e_{\theta} +\frac{1}{\mathrm{sin}\theta}(\partial_{\varphi}\Phi_{s})e_{\varphi}.
\end{eqnarray*}
We define the divergence of $v$ by
\begin{eqnarray*}
\mathrm{div} v= \mathrm{div} (v_{\theta}e_{\theta}+v_{\varphi}e_{\varphi} )=\frac{1}{\mathrm{sin}\theta}(\partial_{\theta}(v_{\theta}\mathrm{sin} \theta)+\partial_{\varphi} v_{\varphi}  ).
\end{eqnarray*}
The horizontal Laplace-Beltrami operator of scalar functions $T$ and $q$ are
\begin{eqnarray*}
\Delta T= \mathrm{div}(\mathrm{grad}T)=\frac{1}{\mathrm{sin}\theta}[\partial_{\theta}(\mathrm{sin}\theta\partial_{\theta}T )+\frac{1}{\mathrm{sin}\theta}\partial_{\varphi\varphi}T],
\end{eqnarray*}
\begin{eqnarray*}
\Delta q= \mathrm{div}(\mathrm{grad}q)=\frac{1}{\mathrm{sin}\theta}[\partial_{\theta}(\mathrm{sin}\theta\partial_{\theta}q )+\frac{1}{\mathrm{sin}\theta}\partial_{\varphi\varphi}q].
\end{eqnarray*}
We define the horizontal Laplace-Beltrami operator $\Delta$ for vector functions on $S^{2}$ as
\begin{eqnarray*}
\Delta v=(\Delta v_{\theta}-\frac{2\mathrm{cos}\theta}{\mathrm{sin}^{2}\theta}\partial_{\varphi}v_{\varphi}-\frac{v_{\theta}}{\mathrm{sin}^{2}\theta  }  )e_{\theta}+(\Delta v_{\varphi}+ \frac{2\mathrm{cos}\theta}{\mathrm{sin}^{2}\theta}\partial_{\varphi}v_{\theta}-\frac{v_{\varphi}}{\mathrm{sin}^{2}\theta   } )e_{\varphi}.
\end{eqnarray*}
Consequently, by integration by parts, we have
\begin{eqnarray}
\int_{0}^{1}w\partial_{\xi}vd\xi= \int_{0}^{1}v \mathrm{div} v d\xi= \int_{0}^{1}\tilde{v} \mathrm{div} \tilde{v}d\xi,
\end{eqnarray}
\begin{eqnarray}
\int_{0}^{1}\nabla_{v}vd\xi= \int_{0}^{1}\nabla_{\tilde{v}}\tilde{v}d\xi+\nabla_{\bar{v}}\bar{v}.
\end{eqnarray}
Taking the average of equations $(1.11)$ in the $z$ direction, over the interval $(0,1)$ and using $(2.18)-(2.20)$ and the boundary conditions $(1.15 )-(1.16),$  we arrive at
\begin{eqnarray}
\partial_{t}\bar{v}&+&\nabla_{\bar{v}}\bar{v}+\overline{\tilde{v}\mathrm{div}\tilde{v} +  \nabla_{\tilde{v}}\tilde{v}}+\frac{f}{R_{0}} \bar{v}^{\bot}+\mathrm{grad } \Phi_{s}+\int_{0}^{1}\int_{\xi}^{1}\frac{bP}{p}\mathrm{grad}[(1+aq)T ]d\xi'd\xi\nonumber\\
&&-\Delta \bar{v}=0\ \ \mathrm{in}\ S^{2}.
\end{eqnarray}
By subtracting $(2.21)$ from $(1.11),$ we obtain the following equation
\begin{eqnarray}
\partial_{t}\tilde{v}&+&\nabla_{\tilde{v}} \tilde{v}+\Big{(}\int_{\xi}^{1} \mathrm{div} \tilde{v} d\xi' \Big{)}\partial_{\xi} \tilde{v}
+\nabla_{\tilde{v}}\bar{v}+\nabla_{\bar{v}}\tilde{v}-\overline{(\tilde{v}  \mathrm{div} \tilde{v} + \nabla_{\tilde{v}}\tilde{v}  )}
+\frac{f}{R_{0}}\tilde{v}^{\bot}\nonumber\\
&&+\int_{\xi}^{1}\frac{bP}{p} \mathrm{grad}[(1+aq)T  ]d\xi'- \int_{0}^{1}\int_{\xi}^{1}\frac{bP}{p} \mathrm{grad}[(1+aq)T  ]d\xi'd\xi\nonumber\\
&&-\Delta \tilde{v}-\partial_{\xi\xi}\tilde{v}=0\ \ \mathrm{in}\ \mho,
\end{eqnarray}
with the following boundary value conditions
\begin{eqnarray}
\partial_{\xi} \tilde{v}=0\ \mathrm{on}\ \xi=1\ \mathrm{and}\  \xi=0.
\end{eqnarray}
Let $e_{\theta}, e_{\varphi}, e_{\xi}$ be the unite vectors in $\theta, \varphi$ and $\xi$ directions of the space domain $\mho$ respectively,
\begin{eqnarray*}
e_{\theta}=\partial_{\theta},\ \ \ \ e_{\varphi}= \frac{1}{\mathrm{sin\theta}}\partial_{\varphi},\ \ \ \ e_{\xi}=\partial_{\xi}.
\end{eqnarray*}
The inner product and norm on $T_{(\theta, \varphi, \xi)}\mho$  (the tangent space of $\mho$ at the point $(\theta,  \varphi, \xi )$  ) are defined by
\begin{eqnarray*}
(u,v)=u\cdot v=\sum_{i=1}^{3} u_{i}v_{i},\ \ \ \ |u|=(u,\ u)^{\frac{1}{2}},
\end{eqnarray*}
where $u=u_{1}e_{\theta}+u_{2}e_{\varphi}+u_{3}e_{\xi}\in T_{(\theta, \varphi, \xi)}\mho$ and  $v=v_{1}e_{\theta}+v_{2}e_{\varphi}+v_{3}e_{\xi}\in T_{(\theta, \varphi, \xi)}\mho.$ For $1\leq p\leq \infty,$ let $L^{p}(\mho), L^{p}(S^{2})$ be the usual Lebesgue spaces with the norm $|\cdot|_{p}$ and $|\cdot|_{L^{p}(S^{2})}$
respectively. If there is no confusion, we will write $|\cdot|_{p}$ instead of $|\cdot|_{L^{p}(S^{2})}$.  $L^{2}(T\Omega|TS^{2} ) $ is the first two components of $L^{2}$ vector fields on $\mho$ with the norm $|v|_{2}=(\int_{\mho}(|v_{\theta}|^{2}+|v_{\varphi}|^{2} ) d\mho )^{\frac{1}{2}},$  where $v=(v_{\theta}, v_{\varphi} ): \mho\rightarrow TS^{2}$. Denoted by $C^{\infty}(S^{2})$ the functions of all smooth functions from $S^{2}$ to $\mathbb{R}.$ Similarly, we can define $C^{\infty}(\mho)$. $H^{m}(\mho)$ is the Sobolev space of functions which are in $L^{2},$ together with all their covariant derivatives with respect to $e_{\theta}, e_{\varphi}, e_{\xi}$ of order $\leq m,$ with the norm
\begin{eqnarray*}
\|h\|_{m}=[\int_{\mho}(\sum_{1\leq k\leq m}\sum_{i_{j}=1,2,3;j=1,...,k}|\nabla_{i_{1}}\cdots \nabla_{i_{k}}h |^{2}+|h|^{2}   )  ]^{\frac{1}{2}},
\end{eqnarray*}
where $\nabla_{1}=\nabla_{e_{\theta}}, \nabla_{2}=\nabla_{e_{\varphi}}$ and $\nabla_{3}=\partial_{\xi} $ which are defined above. Denote $H^{m}(T\mho|TS^{2}) =\{v;v=(v_{\theta}, v_{\varphi}):\mho\rightarrow TS^{2}, \|v\|_{m}< \infty   \},$ where the norm  is similar to that of $H^{m}(\mho)$ (i.e., let $h=(v_{\theta}, v_{\varphi} )= v_{\theta}e_{\theta}+v_{\varphi}e_{\varphi}.$)

We will conduct our work in the following functional spaces. Let
\begin{eqnarray*}
\mathcal{V}_{1}:=\{v; v\in C^{\infty}(T\mho|TS^{2}),\ \partial_{\xi}v|_{\xi=0}=0,\ \partial_{\xi}v|_{\xi=1}=0,\ \int_{0}^{1}\mathrm{div}\ vd\xi=0  \},
\end{eqnarray*}
\begin{eqnarray*}
\mathcal{V}_{2}:=\{T; T\in C^{\infty}(\mho),\ \partial_{\xi}T|_{\xi=0}=0,\ \partial_{\xi}T|_{\xi=1}=-\alpha_{s}T  \},
\end{eqnarray*}
\begin{eqnarray*}
\mathcal{V}_{3}:=\{q; q\in C^{\infty}(\mho),\ \partial_{\xi}q|_{\xi=0}=0,\ \partial_{\xi}q|_{\xi=1}=-\beta_{s}q  \}.
\end{eqnarray*}
We denote by $V_{1}, V_{2}$ and $V_{3}$ the closure spaces of $\mathcal{V}_{1},  \mathcal{V}_{2}$ and $\mathcal{V}_{3}$ in $H^{1}(\mho)$ under $H^{1}-$ topology, respectively. In addition, we denote by $H_{1}, H_{2}$ and $ H_{3}$ the closure of $\mathcal{V}_{1}, \mathcal{V}_{2}   $  and $\mathcal{V}_{3} $  in $L^{2}(\mho)$  under $L^{2}-$ topology. Let $H:= H_{1}\times H_{2}\times H_{3}$  and $V=V_{1}\times V_{2}\times V_{3}$ with $V'$ being  dual space of $V$.
 By definition, the inner products and norms on $V_{1}, V_{2}$ and $ V_{3}$ are given by
\begin{eqnarray*}
\langle v,v_{1}      \rangle_{V_1}=\int_{\mho}(\nabla_{e_{\theta}}v\cdot \nabla_{e_{\theta}}v_{1}+\nabla_{e_{\varphi}}v\cdot \nabla_{e_{\varphi}}v_{1}+\partial_{\xi}v\partial_{\xi}v_{1}+v\cdot v_{1} )d\mho,
\end{eqnarray*}
\begin{eqnarray*}
\|v\|_{1}=\langle v,     v\rangle_{V_{1}}^{\frac{1}{2}},\ \ \ \ \forall\ v,\ v_{1}\in V_{1},
\end{eqnarray*}
\begin{eqnarray*}
\langle T,   T_{1} \rangle_{V_{2}}=\int_{\mho}(\mathrm{grad} T\cdot \mathrm{grad} T_{1} +\partial_{\xi}T \partial_{\xi}T_{1} )d\mho +\alpha_{s}\int_{S^{2}}TT_{1} dS^{2},
\end{eqnarray*}
\begin{eqnarray*}
\|T\|_{1}=\langle T,   T \rangle_{V_{2}}^{\frac{1}{2}},\ \ \ \forall\ T,T_{1}\in V_{2},
\end{eqnarray*}
\begin{eqnarray*}
\langle q,   q_{1} \rangle_{V_{3}}=\int_{\mho}(\mathrm{grad} q\cdot \mathrm{grad} q_{1} +\partial_{\xi}q \partial_{\xi}q_{1})d\mho + \beta_{s}\int_{S^{2}}qq_{1}dS^{2} ,
\end{eqnarray*}
\begin{eqnarray*}
\|q\|_{1}=\langle q,   q\rangle_{V_{3}}^{\frac{1}{2}},\ \ \ \forall\ q,q_{1}\in V_{3}.
\end{eqnarray*}
  Let $V_{i}'(i=1,2,3)$ be the dual space of $V_{i}$ with $\langle , \rangle$ being the inner products between $V_{i}'$ and $V_{i}.$ Without confusion, we also  denote by $\langle, \rangle$ the inner product in $L^{2}(\mho)$ and $L^{2}(S^{2}).$  Define the linear operator $A_{i}: V_{i}\mapsto V_{i}', i=1,2,3$ :
\begin{eqnarray*}
&&\langle A_{1}u_{1}, u_{2}    \rangle= \langle u,  v  \rangle_{V_{1}},\ \ \ \forall\ u_{1},\ u_{2}\in V_{1};\\
&&\langle A_{2}\theta_{1} ,  \theta_{2}   \rangle= \langle \theta_{1},  \theta_{2}   \rangle_{V_{2}},\ \ \ \forall\ \theta_{1},\theta_{2} \in V_{2};\\
&&\langle A_{3} q_{1} ,  q_{2}   \rangle= \langle q_{1},  q_{2}   \rangle_{V_{3}},\ \ \ \forall\ q_{1},q_{2} \in V_{3}.
\end{eqnarray*}
Denote $D(A_{i})=\{\eta\in V_{i},    A_{i}\eta\in H_{i} \}.$ Since $ A_{i}$ is positive self-adjoint with compact resolvent, according to the classic spectral theory we can define the power $A_{i}^{s}$ for any $s\in \mathbb{R}.$ Then we have $D(A_{i}^{\frac{1}{2}})=V_{i}$ and $D(A_{i}^{-\frac{1}{2}})=V_{i}'.$ Moreover,
\begin{eqnarray*}
D(A_{i})\subset V_{i}\subset H_{i}\subset V_{i}'\subset D(A_{i})',
\end{eqnarray*}
where $D(A_{i})' $ is the dual space of $ D(A_{i})$ and the embeddings above are all compact.
In the following, we state some lemmas including integrations by parts and uniform Gronwall lemma, which are frequently used in our paper. For the proof of Lemma $2.1$-Lemma $2.3$, we can see $\cite {GH1}.$ The proof of uniform Gronwall lemma was given in $\cite{FP, T}.$
\begin{lemma}
Let $u=(u_{\theta}, u_{\varphi} ), v= (v_{\theta}, v_{\varphi})\in C^{\infty}(T\mho|TS^{2} )$ and $p\in C^{\infty}(S^{2}).$ Then
\begin{eqnarray*}
\int_{S^{2}}p\ \mathrm{div}\ udS^{2}=-\int_{S^{2}}\nabla p \cdot u dS^{2},
\end{eqnarray*}
\begin{eqnarray*}
\int_{\mho}\nabla p \cdot v d\mho=0\ \ \ for\ any\ v\in V_{1},
\end{eqnarray*}
and
\begin{eqnarray*}
\int_{\mho}(-\Delta u) \cdot vd\mho= \int_{\mho}(\nabla_{e_{\theta} } u\cdot \nabla_{e_{\theta}}v+ \nabla_{e_{\varphi} } u\cdot \nabla_{e_{\varphi}}v+u\cdot v    )d\mho.
\end{eqnarray*}
\end{lemma}
\begin{lemma}
For any $h\in C^{\infty}(S^{2}), v\in C^{\infty}(T\mho|TS^{2}),$ we have
\begin{eqnarray*}
\int_{S^{2}}\nabla_{v}hd S^{2}+\int_{S^{2}}h\ \mathrm{div}\ v dS^{2} = \int_{S^{2}}\mathrm{div} (hv)dS^{2}=0.
\end{eqnarray*}
\end{lemma}
\begin{lemma}
Let $u,v \in V_{1}, T\in V_{2}, q\in V_{3}.$ Then we have
\begin{eqnarray*}
\int_{\mho}[\nabla_{u}v+(\int_{\xi}^{1} \mathrm{div }\ u d\xi'  )\partial_{\xi}v ]v d\mho=0,
\end{eqnarray*}
\begin{eqnarray*}
\int_{\mho}[\nabla_{u}g+(\int_{\xi}^{1} \mathrm{div}\ u d\xi'  )\partial_{\xi}g ]g d\mho=0,\ \ for\ g=T\ or\ g=q,
\end{eqnarray*}
\begin{eqnarray*}
\int_{\mho}\Big{(} \int_{\xi}^{1}\frac{bP}{p} \mathrm{grad} [(1+aq)T]d\xi'\cdot u-\frac{bP}{p}(1+aq)T(\int_{\xi}^{1} \mathrm{div }\ u d\xi'  )   \Big{)}=0.
\end{eqnarray*}
\end{lemma}
In our article, we will frequently use the following inequalities. So we state them as the lemmas below. For their proof, one can refer to $\cite{GH2}$ and $\cite{ZG}.$
\begin{lemma}
Let $v\in H^{2}(T\mho|TS^{2}), \mu\in H^{1}(T\mho|TS^{2})  \Big{(}\mu\in H^{1}(\mho)\Big{)}$ and $\nu\in  L^{2}(T\mho|TS^{2})  \Big{(} \nu \in L^{2}(\mho)\Big{)}.$ Then, there exists a positive constant $c$ independent of $v,\mu$ and $\nu$ such that
\begin{eqnarray*}
&&\Big{|}\langle \Big{(} \int_{\xi}^{1} \mathrm{div} v(t; \theta, \phi, \xi' )d\xi'  \Big{)} \mu,  \nu \rangle\Big{|}\\
&\leq& c|\mathrm{div} v|_{2}^{\frac{1}{2}}(|\mathrm{div} v|_{2}^{\frac{1}{2}}+ |\Delta v|_{2}^{\frac{1}{2}} )|\mu|_{2}^{\frac{1}{2}}(|\nabla_{e_{\theta}} \mu|_{2}^{\frac{1}{2}}+|\nabla_{e_{\varphi}} \mu|_{2}^{\frac{1}{2}}   + |\Delta \mu|_{2}^{\frac{1}{2}})|\nu|_{2}\\
\Big{(}&\leq& c|\mathrm{div} v|_{2}^{\frac{1}{2}}(|\mathrm{div} v|_{2}^{\frac{1}{2}}+ |\Delta v|_{2}^{\frac{1}{2}} )|\mu|_{2}^{\frac{1}{2}}(|\nabla \mu|_{2}^{\frac{1}{2}} + |\Delta \mu|_{2}^{\frac{1}{2}})|\nu|_{2}\Big{)}.
\end{eqnarray*}
\end{lemma}
 \begin{lemma}
Let $v\in H^{1}(T\mho|TS^{2}), \mu\in H^{1}(T\mho|TS^{2}) \Big{(}\mu\in H^{1}(\mho)\Big{)}$ and   $\nu\in H^{1}(T\mho|TS^{2})     \Big{(}  \nu\in H^{1}(\mho)\Big{)}.$ Then, there exists a positive constant $c$ independent of $v,\mu$ and $\nu$ such that
\begin{eqnarray*}
&&\Big{|}\langle \Big{(} \int_{\xi}^{1} \mathrm{div} v(t; \theta, \phi, \xi' )d\xi'  \Big{)} \mu,  \nu \rangle\Big{|}\\
&\leq& c| \mathrm{div} v|_{2}|\mu|_{2}^{\frac{1}{2}}( |\mu|_{2}^{\frac{1}{2}}+ |\nabla_{e_{\theta}}\mu|_{2}^{\frac{1}{2}}+ |\nabla_{e_{\varphi}}\mu|_{2}^{\frac{1}{2}} )
|\nu|_{2}^{\frac{1}{2}}( |\nu|_{2}^{\frac{1}{2}}+ |\nabla_{e_{\theta}}\nu|_{2}^{\frac{1}{2}}+ |\nabla_{e_{\varphi}}\nu|_{2}^{\frac{1}{2}} )\\
\Big{(} &\leq& c| \mathrm{div} v|_{2}|\mu|_{2}^{\frac{1}{2}}( |\mu|_{2}^{\frac{1}{2}}+ |\nabla\mu|_{2}^{\frac{1}{2}})
|\nu|_{2}^{\frac{1}{2}}( |\nu|_{2}^{\frac{1}{2}}+ |\nabla\nu|_{2}^{\frac{1}{2}})   \Big{)}.
\end{eqnarray*}
\end{lemma}
In the process of obtaining absorbing ball for the strong solution $U$ in $H^{2}(T\mho|TS^{2})\times H^{2}(\mho)\times H^{2}(\mho),$ the uniform Gronwall lemma is used extensively. Therefore, for sake of convenience,  we cite it here.
\begin{lemma}
 Let $f,g$ and $h$ be three non-negative locally integrable functions on $(t_{0}, \infty)$ such that
\begin{eqnarray*}
 \frac{df}{dt}\leq gf +h,\ \ \ \forall\ t\geq t_{0},
\end{eqnarray*}
and
\begin{eqnarray*}
\int_{t}^{t+r}f(s)ds\leq a_{1},\ \ \int_{t}^{t+r}g(s)ds\leq a_{2},\ \ \int_{t}^{t+r}h(s)ds\leq a_{3},\ \ \forall\ t\geq t_{0},
\end{eqnarray*}
where $r, a_{1}, a_{2}, a_{3}$ are positive constants. Then
\begin{eqnarray*}
 f(t+r)\leq (\frac{a_{1}}{r}+a_{3} )e^{a_{2}},\ \ \ \forall\ t\geq t_{0}.
\end{eqnarray*}
\end{lemma}
Before considering the dimensions of the global attractor for the strong solutions to the moist primitive equations, we recall the definitions of strong solution to $(1.11)-(1.17).$
\begin{definition}
Suppose $Q_{1},Q_{2}\in L^{2}(\mho)$ and $ Q_{1}|_{\xi=1}, Q_{2}|_{\xi=1}\in L^{2}(S^{2}).$ Let $ U_{0}:=(v_{0}, T_{0}, q_{0})\in V$ and $\tau>0.$ $U:=(v,T, q)$ is called a strong solution of $(1.11)-(1.17)$ on the time interval $[0, \tau]$ if it satisfies $(1.11)-(1.13)$ in a week sense, and also
\begin{eqnarray*}
v\in C([0, \tau]; V_{1})\cap L^{2}([0,\tau]; H^{2}(\mho)),
\end{eqnarray*}
\begin{eqnarray*}
T\in C([0, \tau]; V_{2})\cap L^{2}([0,\tau]; H^{2}(\mho)),
\end{eqnarray*}
\begin{eqnarray*}
q\in C([0, \tau]; V_{3})\cap L^{2}([0,\tau]; H^{2}(\mho)),
\end{eqnarray*}
\begin{eqnarray*}
\partial_{t} v, \partial_{t} T,  \partial_{t} q \in L^{1}([0, \tau]; L^{2}(\mho)).
\end{eqnarray*}
\end{definition}
Now we state the global well-posedness theorem for the strong solution as follows. For the proof of the therem, one can refer to $\cite{GH2}.$
\begin{proposition}
Let $Q_{1}, Q_{2} \in H^{1}(\mho), U_{0}=(v_{0}, T_{0}, q_{0}) \in V .$
Then for any $\tau > 0$ given, the global strong solution $U$ of the system $(1.11)-(1.17)$ is unique on the interval $[0, \tau]$. Moreover, the strong solution $U$ is continuous with respect to initial data in $H.$
\end{proposition}
\begin{Rem}
In fact, we proved in $\cite{ZG}$ that for any $Q_{1}, Q_{2}\in L^{2}(\mho)$ and $ Q_{1}|_{\xi=1}, Q_{2}|_{\xi=1}\in L^{2}(S^{2}),$ the global strong solution $U$ of the system $(1.11)-(1.17)$ is unique on the interval $[0, \tau]$ with any $\tau>0$. Moreover, the strong solution $U$ is continuous with respect to initial data in $V.$
\end{Rem}
For the reader's convenience,   we introduce the definition of global attractor in the following. For more details, we refer to $\cite{H, T}$ and other references.
Let $(X,d)$ be a separable metric space and $ S(t): X\rightarrow X, 0\leq t<\infty,$ be a semigroup satisfying:
\par
$(i)\ S(t)S(s)x=S(t+s)x, $ for all $t,s \in \mathbb{R}_+$ and $x\in X;$
\par
$(ii)\ S(0)=I$ (Identity in $X$);
\par
$(iii)\ S(t)$ is continuous in $X$ for all $t\geq0 .$
\par
Typically, $S(t)$ is associated with a autonomous differential equation; $S(t)x$ is the state at time $t$ of the solution whose initial data is $x.$

\begin{definition}
   A subset $\mathcal{A}$ in $X$ is said to be a global attractor if it satisfies the following properties:\\
$(i)$ $\mathcal{A}$ is compact in $X;$\\
$(ii)$ for every $t\geq 0, S(t) \mathcal{A} =\mathcal{A};$\\
$(iii)$ for every bounded set $B$ in $X,$ the set $S(t)B$ converges to $\mathcal{A}$ in $X$, when $t\rightarrow \infty, i.e.,$
$$\lim\limits_{t\rightarrow \infty}d(S(t)B, \mathcal{A}) =0. $$
Here, and in the following, for $A$and $ B$ subsets of $X, d(A, B)$ is the semi-distance given by
$$d(A,B)={\sup\limits_{x\in A}}\inf\limits_{y\in B} d(x,y).$$
\end{definition}
In $\cite{ZG},$  we prove the result below about the existence of global attractor for the solutions to 3D moist primitive equations.
\begin{theorem}
Assume $Q_{1}, Q_{2}\in L^{2}(\mho)$ and $ Q_{1}|_{\xi=1}, Q_{2}|_{\xi=1}\in L^{2}(S^{2}).$ Then, for $t\geq 0,$ the solution operator $\{S(t)   \}_{t\geq 0}$ of the $3D$ viscous PEs of large-scale moist atmosphere $(1.11)-(1.17): S(t)(v_{0}, T_{0}, q_{0})= (v(t), T(t), q(t) )$ defines a semigroup in the space $V.$ Furthermore, the results below hold:
\par
$(1)$ For any $(v_{0}, T_{0}, q_{0})\in V, t\rightarrow S(t)(v_{0}, T_{0}, q_{0}) $ is a continuous map from $\mathbb{R}_{+} $ into $V.$
\par
$(2)$ For any $t>0, S(t)$ is a continuous map in $V.$
\par
$(3)$ For any $t>0, S(t)$ is a compact map in $V.$
\par
$(4)$ $\{S(t)\}_{t\geq 0}$ possesses a global attractor $\mathcal{A}$ in $V.$ The global attractor in $\mathcal{A}$ is compact and

\quad connected in $V$ and is the maximal bounded attractor in $V$ in the sense of set inclusion

\quad  relation; $\mathcal{A }$ attracts all bounded subset in $V$ in the norm of $V.$
\end{theorem}
Therefore, the main result of this work is the finite Hausdorff and fractal dimensions of the global attractor in space $V.$ To prove the properties of global attractor,  we use a theorem in $\cite{L0}.$ For reader's convenience, we cite it below.
\begin{theorem}
Let $X$ be a Hilbert space with norm $\|\cdot\|_{X}$. $S:X\mapsto X$ be a map and $\mathcal{A}\subset X$ be a compact set such that $S(\mathcal{A})=\mathcal{A} .$ Suppose that there exist a positive constant $c$ and $\delta \in (0, 1),$ such that $\forall a_{1}, a_{2}\in \mathcal{A}.$
\begin{eqnarray*}
&\mathrm{(i)}&\ \|S(a_{1})- S(a_{2})\|_{X}\leq c\|a_{1}-a_{2}\|_{X},\\
&\mathrm{(ii)}&\ \|Q_{N}[S(a_{1})- S(a_{2})]\|_{X}\leq \delta\|a_{1}-a_{2}\|_{X},
\end{eqnarray*}
where $Q_{N}$ is the projection in $X$ onto some subspace $(X_{N})^{\bot}$ of co-dimension $N\in \mathbb{N}.$ Then
\begin{eqnarray*}
d_{H}(\mathcal{A}) \leq d_{F}(\mathcal{A})\leq N\frac{\mathrm{ln}(\frac{8{G_{a}}^{2}c^{2}}{1-\delta^{2}})}{\mathrm{ln}(\frac{2}{1+\delta^{2}} )},
\end{eqnarray*}
where $d_{H}(\mathcal{A})$ and $d_{F}(\mathcal{A}) $ are the Hausdorff and fractal dimensions of $\mathcal{A}$ respectively and $G_{a}$ is the Gauss constant:
\begin{eqnarray*}
G_{a}:=\frac{1}{2\pi}\beta(\frac{1}{4}, \frac{1}{2} )=\frac{2}{\pi}\int_{0}^{1}\frac{dx}{\sqrt{1-x^4}}=0.8346268....
\end{eqnarray*}

\section{Uniform estimates and absorbing balls for $\partial_{t}U$ in the space $H$}
\begin{theorem}
Assume $Q_{1}, Q_{2}\in L^{2}(\mho)$ and $ Q_{1}|_{\xi=1}, Q_{2}|_{\xi=1}\in L^{2}(S^{2}).$ Let $U_{0}:=(v(0), T(0), q(0))\in V$ and $\partial_{t}U(0)\in H.$ Then, there exists a unique strong solution $U:=(v,T,q)$ of $(1.11)-(1.17)$ such that
\begin{eqnarray*}
\partial_{t}U \in L^{\infty}([0, \infty); H ) \cap L^{2}([0, \infty); V ).
\end{eqnarray*}
Furthermore, there exists a bounded absorbing ball for $U_{t}:=\partial_{t}U$ in space $H. $
\end{theorem}
\begin{proof}
The uniqueness of such solutions follows by Proposition 2.1. To obtain uniform estimates of $\partial_{t}U$ in $H,$  we
introduce some notations denoted by
\begin{eqnarray*}
u:=v_{t}=\partial_{t}v,\  \theta:=T_{t}=\partial_{t}T,\  \eta:=q_{t}=\partial_{t}q.
\end{eqnarray*}
To prove $\partial_{t}U$ is uniformly bounded in $H$ with respect to $t,$ we first prove that there exists a positive constant c independent of $t$ such that
\begin{eqnarray}
\int_{t}^{t+1}|\partial_{s}U(s)|_{2}^{2}ds<c
\end{eqnarray}
for arbitrary $t\geq 0.$
By Lemma 2.1, we have
\begin{eqnarray}
\int_{\mho} q_{t}\Delta qd\mho= -\int_{\mho} \nabla q_{t}\cdot \nabla qd\mho=-\frac{1}{2}\frac{d}{dt}|\nabla q|_{2}^{2}.
\end{eqnarray}
By the boundary conditions $(1.15)-(1.16)$, we arrive at
\begin{eqnarray}
\int_{\mho}q_{t}\partial_{\xi\xi}qd\mho&=&\int_{\mho}[\partial_{\xi} (q_{t} q_{\xi} )-(\partial_{\xi}q_{t})q_{\xi}]d\mho\nonumber\\
&=&\int_{S^{2}}-\beta_{s}q_{t}|_{\xi=1}q|_{\xi=1}dS^{2}-\frac{1}{2}\frac{d}{dt}|q_{\xi}|_{2}^{2}\nonumber\\
&=&-\frac{\beta_{s}}{2}\frac{d}{dt}|q|_{\xi=1}|_{2}^{2}-\frac{1}{2}\frac{d}{dt}|q_{\xi}|_{2}^{2}.
\end{eqnarray}
Using $(3.25)-(3.26)$ and taking the inner product of equation $(1.13)$ with $q_{t}$ yields
\begin{eqnarray}
\frac{1}{2}\frac{d}{dt}|\nabla q|_{2}^{2}&+&\frac{1}{2}\frac{d}{dt}| q_{\xi}|_{2}^{2}+\frac{\beta_{s}}{2}\frac{d}{dt}| q|_{\xi=1}|_{2}^{2}+|\eta|_{2}^{2}\nonumber\\
&=&-\langle \nabla_{v}q,  \eta\rangle-\langle \Big{(} \int_{\xi}^{1}\mathrm{div}\ v( t;\theta,\varphi,\xi')d\xi'\Big{)}\partial_{\xi}q , \eta   \rangle+ \langle Q_{2}, \eta\rangle\nonumber\\
&\leq &|\eta|_{2}|v|_{4}|\nabla q|_{4}+c|\eta|_{2}\|v\|_{1}^{\frac{1}{2}}\|v\|_{2}^{\frac{1}{2}}\|q\|_{1}^{\frac{1}{2}}\|q\|_{2}^{\frac{1}{2}}+|Q_{2}|_{2}|\eta|_{2}\nonumber\\
&\leq&\varepsilon |\eta|_{2}^{2}+c|Q_{2}|_{2}^{2}+c\|v\|_{1}^{2}\|q\|_{2}^{2}+c\|v\|_{1}\|v\|_{2}\|q\|_{1}\|q\|_{2},
\end{eqnarray}
where the first inequality follows by Lemma $2.4.$ Then $(3.27)$ implies that for $t>0$
\begin{eqnarray}
\|q(t+1)\|_{1}^{2}&+&\int_{t}^{t+1}|\eta(s)|_{2}^{2}ds\leq \|q(t)\|_{1}^{2}+|Q_{2}|_{2}^{2}\nonumber\\
&&+c\int_{t}^{t+1}(\|v(s)\|_{1}^{2}\|q(s)\|_{2}^{2}+c\|v(s)\|_{1}\|v(s)\|_{2}\|q(s)\|_{1}\|q(s)\|_{2})ds.
\end{eqnarray}
It is proved in $\cite{ZG}$ that there exists a positive constant $c$ independent of $t$ such that
\begin{eqnarray}
\|U(t)\|_{1}^{2}+\int_{t}^{t+1}\|U(s)\|_{2}^{2}ds<c
\end{eqnarray}
for arbitrary $t\geq0.$
Therefore, by $(3.28)$ and $(3.29)$,  we obtain the uniform boundedness for $ \int_{t}^{t+1}|\eta(s)|_{2}^{2}ds$ with respect $t.$ Taking an similar argument of $(3.27),$ we have
\begin{eqnarray}
\frac{1}{2}\frac{d}{dt}|\nabla T|_{2}^{2}&+&\frac{1}{2}\frac{d}{dt}| T_{\xi}|_{2}^{2}+\frac{\alpha_{s}}{2}\frac{d}{dt}| T|_{\xi=1}|_{2}^{2}+|\theta|_{2}^{2}\nonumber\\
&=&-\langle \nabla_{v}T,  \theta\rangle-\langle \Big{(} \int_{\xi}^{1}\mathrm{div}\ v( t;\theta,\varphi,\xi')d\xi'\Big{)}\partial_{\xi}T , \theta   \rangle\nonumber\\
&&+\langle Q_{1},    \theta  \rangle+ \langle \frac{bP}{p}(1+aq)  \Big{(} \int_{\xi}^{1}\mathrm{div}\ v( t;\theta,\varphi,\xi')d\xi'\Big{)}  ,        \theta  \rangle\nonumber\\
&\leq& |\theta|_{2}|\nabla T|_{4}|v|_{4}+c|\theta|_{2}\|T\|_{1}^{\frac{1}{2}}\|T\|_{2}^{\frac{1}{2}}\|v\|_{1}^{\frac{1}{2}}\|v\|_{2}^{\frac{1}{2}}\nonumber\\
&&+|\theta|_{2}|Q_{1}|_{2}+c|\theta|_{2}\|v\|_{1}+c|\theta|_{2}\|q\|_{1}^{\frac{1}{2}}\|q\|_{2}^{\frac{1}{2}}\|v\|_{1}^{\frac{1}{2}}\|v\|_{2}^{\frac{1}{2}}\nonumber\\
&\leq&\varepsilon |\theta|_{2}^{2}+c\|T\|_{2}^{2}\|v\|_{1}^{2}+c\|T\|_{1}\|T\|_{2}\|v\|_{1}\|v\|_{2} \nonumber\\
&&+c|Q_{1}|_{2}^{2}+c\|v\|_{1}^{2}+c\|q\|_{1}\|q\|_{2}\|v\|_{1}\|v\|_{2},
\end{eqnarray}
where the first inequality follows by Lemma $2.4.$ For $t>0,$ integrating $(3.30)$ with respect to $t$ yields
\begin{eqnarray}
&&\|T(t+1)\|_{1}^{2}+\int_{t}^{t+1}|\theta(s)|_{2}^{2}ds\nonumber\\
&\leq&\|T(t)\|_{1}^{2}+c\int_{t}^{t+1}\|T(s)\|_{2}^{2}\|v(s)\|_{1}^{2}ds\nonumber\\
&&+c\int_{t}^{t+1}\|T(s)\|_{1}\|T(s)\|_{2}\|v(s)\|_{1}\|v(s)\|_{2}ds\nonumber\\
&&+c|Q_{1}|_{2}^{2}+c\int_{t}^{t+1}\|v(s)\|_{1}^{2}ds\nonumber\\
&&+c\int_{t}^{t+1}\|q(s)\|_{1}\|q(s)\|_{2}\|v(s)\|_{1}\|v(s)\|_{2}ds,
\end{eqnarray}
which combined with $(3.29)$ implies the uniform boundednss of $ \int_{t}^{t+1}|\theta(s)|_{2}^{2}ds$  with respect to $t$.
Since by Lemma 2.1 and boundary conditions $(1.15)-(1.16)$, we have
\begin{eqnarray}
\int_{\mho}\mathrm{grad} \Phi_{s}(\theta, \phi,t)\cdot ud\mho=\int_{\mho}\mathrm{grad} \Phi_{s}(\theta, \phi,t)\cdot v_{t}d\mho=0,
\end{eqnarray}
\begin{eqnarray}
-\int_{\mho} v_{t}\cdot\Delta v d\mho=\int_{\mho}(\nabla_{e_{\theta}}v_{t} \cdot \nabla_{e_{\theta}}v +\nabla_{e_{\phi}}v_{t}\cdot \nabla_{e_{\phi}}v+v_{t}\cdot v)dt,
\end{eqnarray}
\begin{eqnarray}
-\int_{\mho} v_{t}\cdot\partial_{\xi\xi} v d\mho=\frac{ 1}{2}\frac{d}{dt}|\partial_{\xi}v|_{2}^{2}.
\end{eqnarray}
Then taking inner product of  $(1.11)$ with $v_{t}$ and using equalities $( 3.32)-(3.34)$ yields
\begin{eqnarray}
&&\frac{1}{2}\frac{d}{dt}\|v\|_{1}^{2}+|u|_{2}^{2}\nonumber\\
&=&-\langle \nabla_{v}v,  u \rangle-\langle \Big{(} \int_{\xi}^{1}\mathrm{div}\ v( t;\theta,\varphi,\xi')d\xi'\Big{)}  \partial_{\xi}v,    u \rangle\nonumber\\
&&- \langle\frac{f}{R_{0}}v^{\bot}, u\rangle- \langle \Big{(} \int_{\xi}^{1}\frac{bP}{p}\mathrm{grad}[(1+aq)T]d\xi'\Big{)}, u\rangle\nonumber\\
&\leq& |v|_{\infty}\|v\|_{1}|u|_{2}+c\|v\|_{1}\|v\|_{2}|u|_{2}+|v|_{2}|u|_{2}\nonumber\\
&&+c|\nabla T|_{2}|u|_{2}+c|T|_{4}|\nabla q|_{4}|u|_{2}+c|\nabla T|_{4}|q|_{4}|u|_{2}\nonumber\\
&\leq& \varepsilon |u|_{2}^{2}+c|v|_{2}^{2}+c\|T\|_{1}^{2}+c\|v\|_{1}^{2}\|v\|_{2}^{2}+c\|T\|_{1}^{2}\|q\|_{2}^{2}+c\|q\|_{1}^{2}\|T\|_{2}^{2}.
\end{eqnarray}
Using uniform estimates $(3.29)$ and integrating $(3.25)$ with respect to $t,$ we prove that
$\int_{t}^{t+1} |u(s)|_{2}^{2}ds$
is uniformly bounded with respect to $t( t\geq 0).$ With these time-uniform $a\ priori$ estimates of
\begin{eqnarray*}
\int_{t}^{t+1} |u(s)|_{2}^{2}ds,\ \  \int_{t}^{t+1} |\theta(s)|_{2}^{2}ds\ \ \mathrm{and}\ \int_{t}^{t+1} |\eta(s)|_{2}^{2}ds ,
\end{eqnarray*}
 we prove $(3.24).$ In the following, we will show the existence of the absorbing ball of $U_{t}$ in $H.$ Taking derivative with respect to $t$ in $(1.13),$ we have
\begin{eqnarray}
\eta_{t}+L_{3}\eta+\nabla_{u}q+\nabla_{v}\eta+w_{t}q_{\xi}+w\eta_{\xi}=0.
\end{eqnarray}
Taking inner product of $(3.36)$ with $\eta$ and using boundary conditions $(1.15)-(1.16),$ we have
\begin{eqnarray*}
&&\frac{1}{2}\frac{d}{dt} |\eta|_{2}^{2}+|\nabla \eta|_{2}^{2}+|\eta_{\xi}|_{2}^{2}+\beta_{s}|\eta|_{\xi=1}|_{2}^{2}\\
&=& -\langle \nabla_{u}q, \eta \rangle  -\langle \nabla_{v}\eta,   \eta   \rangle - \langle w_{t}q_{\xi},  \eta  \rangle-\langle w\eta_{\xi},  \eta \rangle\\
&=&-\langle \nabla_{u}q, \eta \rangle- \langle w_{t}q_{\xi},  \eta  \rangle,
\end{eqnarray*}
where the second equality follows by Lemma $2.3. $ Therefore, by H$\mathrm{\ddot{o}}$lder inequality, Lemma $2.4,$ interpolation inequality and Young's inequality we obtain
\begin{eqnarray}
&&\frac{d}{dt}|\eta|_{2}^{2}+\|\eta\|_{1}^{2}\nonumber\\
&\leq& c|\nabla q|_{2}|u|_{4}|\eta|_{4}+c|\mathrm{div} u|_{2}|q_{\xi}|_{2}^{\frac{1}{2}}|\nabla q_{\xi}|_{2}^{\frac{1}{2}}|\eta|_{2}^{\frac{1}{2}}|\nabla\eta|_{2}^{\frac{1}{2}}\nonumber\\
&\leq&c\|q\|_{1}|u|_{2}^{\frac{1}{4}}\|u\|_{1}^{\frac{3}{4}}|\eta|_{2}^{\frac{1}{4}}\|\eta\|_{1}^{\frac{3}{4}}++c|\mathrm{div} u|_{2}|q_{\xi}|_{2}^{\frac{1}{2}}|\nabla q_{\xi}|_{2}^{\frac{1}{2}}|\eta|_{2}^{\frac{1}{2}}|\nabla\eta|_{2}^{\frac{1}{2}}\nonumber\\
&\leq& \varepsilon \|\eta\|_{1}^{2}+\varepsilon \|u\|_{1}^{2}+c|\eta|_{2}^{2}+c\|q\|_{1}^{8}|u|_{2}^{2}+c\|q\|_{1}^{2}\|q\|_{2}^{2}|\eta|_{2}^{2}.
\end{eqnarray}
Taking derivative with respect to $t$ in $(1.12),$ we have
\begin{eqnarray}
\partial_{t}\theta&+&L_{2}\theta+\nabla_{u}T+\nabla_{v}\theta+w_{t}T_{\xi}+w\theta_{\xi}\nonumber\\
&&-\frac{bP}{p}a\eta w-\frac{bP}{p}(1+aq)w_{t}=0.
\end{eqnarray}
Since by Lemma $2.3,$
\begin{eqnarray*}
\langle \nabla_{v}\theta+w\theta_{\xi},   \theta   \rangle=0.
\end{eqnarray*}
Consequently, multiplying $(3.38)$ by $\theta$ and integrating over $\mho$ yields
\begin{eqnarray}
&&\frac{1}{2}\frac{d}{dt}|\theta|_{2}^{2}+|\nabla \theta|_{2}^{2}+|\theta_{\xi}|_{2}^{2}+\alpha_{s}|\theta|_{\xi=1}|_{2}^{2}\nonumber\\
&=&\langle \nabla_{u}T+w_{t}T_{\xi},  \theta    \rangle\nonumber\\
&&-\langle \frac{bP}{p}a\eta w,   \theta    \rangle-\langle \frac{bP}{p}(1+aq)w_{t},     \theta \rangle\nonumber\\
&=:& I_{1}+I_{2}+I_{3}.
\end{eqnarray}
By H$\mathrm{\ddot{o}}$lder inequality, Lemma $2.4$, interpolation inequality and Young's inequality we have
\begin{eqnarray*}
I_{1}&\leq& |u|_{4}|\nabla T|_{4}|\theta|_{2}+c\|u\|_{1}|T_{\xi}|_{2}^{\frac{1}{2}}\|T_{\xi}\|_{1}^{\frac{1}{2}}|\theta|_{2}^{\frac{1}{2}}\|\theta\|_{1}^{\frac{1}{2}}\\
&\leq&c|\theta|_{2}|u|_{2}^{\frac{1}{4}}\|u\|_{1}^{\frac{3}{4}}\|T\|_{1}^{\frac{1}{4}}\|T\|_{2}^{\frac{3}{4}}
+c\|u\|_{1}|T_{\xi}|_{2}^{\frac{1}{2}}\|T_{\xi}\|_{1}^{\frac{1}{2}}|\theta|_{2}^{\frac{1}{2}}\|\theta\|_{1}^{\frac{1}{2}}\\
&\leq& \varepsilon \|u\|_{1}^{2}+\varepsilon\|\theta\|_{1}^{2}+c|u|_{2}^{2}+c|\theta|_{2}^{2}\|T\|_{2}^{2}+c|\theta|_{2}^{2}\|T\|_{1}^{2}\|T\|_{2}^{2}.
\end{eqnarray*}
To estimate $I_{2},$ using H$\mathrm{\ddot{o}}$lder inequality, interpolation inequality and Young's inequality we reach
\begin{eqnarray*}
I_{2}&\leq& c|\theta|_{2}|\eta|_{2}^{\frac{1}{2}}\|\eta\|_{1}^{\frac{1}{2}}\|v\|_{1}^{\frac{1}{2}}\|v\|_{2}^{\frac{1}{2}}\\
& \leq& \varepsilon \|\eta\|_{1}^{2}+c|\eta|_{2}^{2}+c|\theta|_{2}^{2}\|v\|_{1}\|v\|_{2}.
\end{eqnarray*}
Taking an analogous argument as $I_{2},$ we deduce
\begin{eqnarray*}
I_{3}&\leq& c\|u\|_{1}|\theta|_{2}+c|q|_{2}^{\frac{1}{2}}\|q\|_{1}^{\frac{1}{2}}\|u\|_{1}|\theta|_{2}^{\frac{1}{2}}\|\theta\|_{1}^{\frac{1}{2}}\\
&\leq & \varepsilon \|u\|_{1}^{2}+\varepsilon\|\theta\|_{1}^{2} +c|\theta|_{2}^{2}|q|_{2}^{2}\|q\|_{1}^{2}.
\end{eqnarray*}
Substituting the estimates of $I_{1}-I_{3}$ into $(3.39),$ we reach
\begin{eqnarray}
\frac{d}{dt}|\theta|_{2}^{2}+\|\theta\|_{1}^{2}
&\leq&  \varepsilon (\|u\|_{1}^{2}+ \|\theta\|_{1}^{2}+ \|\eta\|_{1}^{2}   )+c(|u|_{2}^{2} + |\eta|_{2}^{2}   )\nonumber\\
&&+c|\theta|_{2}^{2}(\|T\|_{2}^{2}+\|v\|_{2}^{2}+\|q\|_{1}^{4}+\|T\|_{1}^{2}\|T\|_{2}^{2}   ).
\end{eqnarray}
Taking derivative with respect to $t$ in $(1.11),$ we have
\begin{eqnarray}
\frac{d}{dt}u&+&L_{1}u+ \nabla_{u}v+\nabla_{v}u+w_{t}v_{\xi}+wu_{\xi}+\frac{f}{R_{0}}u^{\bot}+\mathrm{grad}(\partial_{t}\Phi_{s} )\nonumber \\
&+&\int_{\xi}^{1}\frac{bP}{p}\mathrm{grad} (a \eta T) d\xi'+\int_{\xi}^{1}\frac{bP}{p}\mathrm{grad}[ (1+aq) \eta ] d\xi'=0.
\end{eqnarray}
By the aid of Lemma $2.3$ and Lemma $2.1,$ it follows that
\begin{eqnarray*}
\langle \nabla_{v}u +wu_{\xi},  u  \rangle= \langle \mathrm{grad}(\partial_{t}\Phi_{s} ), u   \rangle=\langle \frac{f}{R_{0}}u^{\bot}, u \rangle=0.
\end{eqnarray*}
Therefore, multiplying $(3.41 )$ by $u$ and integrating over $\mho$ yields
\begin{eqnarray}
\frac{1}{2}\frac{d}{dt}|u|_{2}^{2}+\|u\|_{1}^{2}&=&-\langle \nabla_{u}v,  u  \rangle-\langle w_{t}v_{\xi},  u  \rangle\nonumber\\
&&-\langle \int_{\xi}^{1}\frac{bP}{p}\mathrm{grad}(a\eta T)d\xi',   u  \rangle-\langle \int_{\xi}^{1}\frac{bP}{p} \mathrm{grad}[(1+aq)\eta ]d\xi', u      \rangle\nonumber\\
&:=& J_{1}+J_{2}+J_{3}+J_{4}.
\end{eqnarray}
In view of H$\mathrm{\ddot{o}}$lder inequality, Sobolev inequality and Young's inequality, we have
\begin{eqnarray*}
J_{1}\leq |u|_{2}|\nabla v|_{4}|u|_{4}\leq \varepsilon\|u\|_{1}^{2}+c|u|_{2}^{2}\|v\|_{2}^{2}.
\end{eqnarray*}
Thanks to Lemma $2.5$ and Young's inequality, we reach
\begin{eqnarray*}
J_{2}&\leq& c\|u\|_{1}\|v\|_{1}^{\frac{1}{2}}\|v\|_{2}^{\frac{1}{2}}|u|_{2}^{\frac{1}{2}}\|u\|_{1}^{\frac{1}{2}}\\
&\leq&\varepsilon \|u\|_{1}^{2}+c|u|_{2}^{2}\|v\|_{1}^{2}\|v\|_{2}^{2}.
\end{eqnarray*}
Similarly, by the Lemma $2.4$, Lemma $2.5$ and Young's inequality, we obtain the estimates of $J_{3}$
\begin{eqnarray*}
J_{3}&\leq& c\|\eta\|_{1}|T|_{2}^{\frac{1}{2}}\|T\|_{1}^{\frac{1}{2}}|u|_{2}^{\frac{1}{2}}\|u\|_{1}^{\frac{1}{2}}+c\|T\|_{1}^{\frac{1}{2}}\|T\|_{2}^{\frac{1}{2}}
|\eta|_{2}^{\frac{1}{2}}\|\eta\|_{1}^{\frac{1}{2}}|u|_{2}\\
&\leq& \varepsilon \|\eta\|_{1}^{2}+\varepsilon \|u\|_{1}^{2}+c\|T\|_{1}^{2}\|T\|_{2}^{2}(|u|_{2}^{2}+|\eta|_{2}^{2}).
\end{eqnarray*}
Finally,  taking an analogously argument, we deduce
\begin{eqnarray*}
J_{4}&\leq& c\|\eta\|_{1}|q|_{2}^{\frac{1}{2}}\|q\|_{1}^{\frac{1}{2}}|u|_{2}^{\frac{1}{2}}\|u\|_{1}^{\frac{1}{2}}+c|u|_{2}\|q\|_{1}^{\frac{1}{2}}\|q\|_{2}^{\frac{1}{2}}
|\eta|_{2}^{\frac{1}{2}}\|\eta\|_{1}^{\frac{1}{2}}+c\|\eta\|_{1}|u|_{2}\\
&\leq&\varepsilon \|\eta\|_{1}^{2}+\varepsilon\|u\|_{1}^{2}+c|u|_{2}^{2}+c\|q\|_{1}^{2}\|q\|_{2}^{2}(|u|_{2}^{2}+|\eta|_{2}^{2}).
\end{eqnarray*}
Substituting these estimates of $J_{1}-J_{4}$ into $(3.42),$ we have
\begin{eqnarray}
&&\frac{d}{dt}|u|_{2}^{2}+\|u\|_{1}^{2}\leq \varepsilon (\|u\|_{1}^{2}+\|\eta\|_{1}^{2} )\nonumber\\
&&+c(|u|_{2}^{2}+|\eta|_{2}^{2})(1+\|v\|_{2}^{2}+\|v\|_{1}^{2}\|v\|_{2}^{2}+
\|T\|_{1}^{2}\|T\|_{2}^{2}+\|q\|_{1}^{2}\|q\|_{2}^{2} ).
\end{eqnarray}
For $t\geq 0,$ denote by
\begin{eqnarray*}
 f(t):= 1+\|q(t)\|_{1}^{8}+\|T(t)\|_{2}^{2}+\|v(t)\|_{2}^{2}+\|q(t)\|_{1}^{2}\|q(t)\|_{2}^{2}+\|T(t)\|_{1}^{2}\|T(t)\|_{2}^{2}+\|v(t)\|_{1}^{2}\|v(t)\|_{2}^{2}.
\end{eqnarray*}
Then, it follows from the time-uniform estimates $(3.29)$ that  $\int_{t}^{t+1}f(s)ds$ is uniformly bounded with respect to $t.$
Summing $(3.37), (3.40)$ and $(3.43),$ we obtain
\begin{eqnarray}
\frac{d}{dt}|\partial_{t}U(t)|_{2}^{2}+\|\partial_{t}U(t)\|_{1}^{2}\leq c|\partial_{t}U(t)|_{2}^{2}f(t).
\end{eqnarray}
Applying the uniform boundedness of $\int_{t}^{t+1}f(s)ds$ and the uniform Gronwall inequality to $(3.44),$ the conclusions of the theorem follows. \hspace{\fill}$\square$
\end{proof}

\section{Dimensions of the global attractor}
\begin{theorem}
Assume $Q_{1}, Q_{2}\in L^{2}(\mho)$ and $ Q_{1}|_{\xi=1}, Q_{2}|_{\xi=1}\in L^{2}(S^{2}).$ Let $U_{0}\in \mathcal{A}.$ Then there is a positive constant $c$ independent of initial condition $U_{0}$ and time $\tau(\tau>0)$ such that
\begin{eqnarray*}
\int_{0}^{\tau}\|U(t)\|_{2}^{2} dt <c(\tau^{\frac{1}{2}} +\tau).
\end{eqnarray*}
\end{theorem}
\begin{proof}
By direct calculation, we reach
\begin{eqnarray*}
||v(\tau)|_{2}^{2}-|v(0)|_{2}^{2}|&=&2|\int_{0}^{\tau}\langle v_{t}(t), v(t)\rangle dt|\leq c\tau.
\end{eqnarray*}
Similarly,
\begin{eqnarray*}
||T(\tau)|_{2}^{2}-|T(0)|_{2}^{2}|&=&2|\int_{0}^{\tau}\langle T_{t}(t),  T(t)\rangle dt|\leq c\tau
\end{eqnarray*}
and
\begin{eqnarray*}
||q(\tau)|_{2}^{2}-|q(0)|_{2}^{2}|&=&2|\int_{0}^{\tau}\langle q_{t}(t), q(t)\rangle dt|\leq c\tau.
\end{eqnarray*}
Taking inner product of $(1.13)$ with $A_{3}q$ and integrating over $[0, \tau]$ with respect to $t$ yields
\begin{eqnarray*}
2\int_{0}^{\tau}\langle q_{t}(t),    A_{3}q(t)   \rangle dt =\|q(\tau)\|_{1}^{2}-\|q(0)\|_{1}^{2}.
\end{eqnarray*}
Consequently, we have
\begin{eqnarray*}
|\|q(\tau)\|_{1}^{2}-\|q(0)\|_{1}^{2}|&\leq& 2\int_{0}^{\tau}|q_{t}(t)|_{2}   | A_{3}q(t)|_{2}  dt \nonumber\\
&\leq &2\sup\limits_{t\in [0, \tau]}|q_{t}(t)|_{2}\int_{0}^{\tau}| A_{3}q(t)|_{2}  dt\nonumber\\
&\leq&c\tau^{\frac{1}{2}}.
\end{eqnarray*}
Similarly, we deduce that
\begin{eqnarray*}
|\|T(\tau)\|_{1}^{2}-\|T(0)\|_{1}^{2}|&\leq& 2\int_{0}^{\tau}|T_{t}(t)|_{2}   | A_{2}T(t)|_{2}  dt \nonumber\\
&\leq &2\sup\limits_{t\in [0, \tau]}|T_{t}(t)|_{2}\int_{0}^{\tau}| A_{2}T(t)|_{2}  dt\nonumber\\
&\leq&c\tau^{\frac{1}{2}}
\end{eqnarray*}
and
\begin{eqnarray*}
|\|v(\tau)\|_{1}^{2}-\|v(0)\|_{1}^{2}|&\leq& 2\int_{0}^{\tau}|v_{t}(t)|_{2}   | A_{1}v(t)|_{2}  dt \nonumber\\
&\leq &2\sup\limits_{t\in [0, \tau]}|v_{t}(t)|_{2}\int_{0}^{\tau}| A_{1}v(t)|_{2}  dt\nonumber\\
&\leq&c\tau^{\frac{1}{2}}.
\end{eqnarray*}
From $\cite{GH2},$ we have
\begin{eqnarray*}
&&\frac{d|v_{\xi}|_{2}^{2}}{dt}+|\nabla_{e_{\theta}}v_{\xi} |_{2}^{2}+|\nabla_{e_{\varphi}}v_{\xi} |_{2}^{2}+|v_{\xi}|_{2}^{2}+|v_{\xi\xi}|_{2}^{2}\nonumber\\
&\leq&c(\|\bar{v}\|_{1}^{8}+|\tilde{v}|_{4}^{8}  )|v_{\xi}|_{2}^{2}+c|T|_{2}^{2}+c|q|_{4}^{4}+c|T|_{4}^{4},
\end{eqnarray*}
which combined with the uniform-time estimates proved in $\cite{ZG}$ implies
\begin{eqnarray}
\int_{0}^{\tau}\|v_{\xi}(t)\|_{1}^{2}dt&\leq& | |v_{\xi}(\tau)|_{2}^{2}- |v_{\xi}(0)|_{2}^{2}  |\nonumber\\
&&+c\int_{0}^{\tau}[(\|\bar{v}(t)\|_{1}^{8}+|\tilde{v}(t)|_{4}^{8}  )|v_{\xi}(t)|_{2}^{2}+|T(t)|_{2}^{2}+|q(t)|_{4}^{4}+|T(t)|_{4}^{4}]dt\nonumber\\
&\leq&c\tau^{\frac{1}{2}}+c\tau.
\end{eqnarray}
In view of the energy estimates proved in $\cite{ZG},$ it follows that
\begin{eqnarray*}
&&\frac{1}{2}\frac{d(|\nabla_{e_{\theta}} v|_{2}^{2}+ |\nabla_{e_{\varphi}} v|_{2}^{2}+|v|_{2}^{2}   )}{dt}+|\Delta v|_{2}^{2}\\
&\leq&c(| \nabla_{e_{\theta}}v |_{2}^{2}+ | \nabla_{e_{\varphi}}v |_{2}^{2}  )(1+|v_{\xi}|_{2}^{4}+|v_{\xi}|_{2}^{2}|\nabla_{e_{\theta}} v_{\xi}|_{2}^{2}+
|v_{\xi}|_{2}^{2}|\nabla_{e_{\varphi}} v_{\xi}|_{2}^{2})\\
&&+c+c|\nabla T|_{2}^{2}+c|\nabla q|_{2}^{2},
\end{eqnarray*}
which implies that
\begin{eqnarray}
\int_{0}^{\tau}|\Delta v(t)|_{2}^{2}dt
&\leq& |\|v(\tau)\|_{1}^{2}- \|v(0)\|_{1}^{2} |+||v(\tau) |_{2}^{2}-|v(0)|_{2}^{2}  |\nonumber\\
&&+ c\int_{0}^{\tau}(| \nabla_{e_{\theta}}v(t) |_{2}^{2}+ | \nabla_{e_{\varphi}}v(t) |_{2}^{2}  )\nonumber\\
&&\cdot(1+|v_{\xi}(t)|_{2}^{4}+|v_{\xi}(t)|_{2}^{2}|\nabla_{e_{\theta}} v_{\xi}(t)|_{2}^{2}+ +|v_{\xi}(t)|_{2}^{2}|\nabla_{e_{\varphi}} v_{\xi}(t)|_{2}^{2} ) dt\nonumber\\
&&+c\int_{0}^{\tau}(1+|\nabla T(t)|_{2}^{2}+|\nabla q(t)|_{2}^{2})dt\nonumber\\
&\leq& c\tau^{\frac{1}{2}}+c\tau.
\end{eqnarray}
According to $\cite{ZG}, $ we know the following estimates
\begin{eqnarray*}
&&\frac{d(|T_{\xi}|_{2}^{2}+|q_{\xi}|_{2}^{2}+\alpha_{s}|T|_{\xi=1}|_{2}^{2}+\beta_{s}|q|_{\xi=1}|_{2}^{2} )}{dt}+ |\nabla T_{\xi}|_{2}^{2}+|T_{\xi\xi}|_{2}^{2}\\
&&+ |\nabla q_{\xi}|_{2}^{2}+|q_{\xi\xi}|_{2}^{2}+\alpha_{s}|\nabla T|_{\xi=1}|_{2}^{2}+\beta_{s}|\nabla q|_{\xi=1}|_{2}^{2}\\
&&\leq c+c(|T_{\xi}|_{2}^{2}+ |q_{\xi}|_{2}^{2}   )+c\|v_{\xi}\|_{1}^{2}+c\|v\|_{1}^{2}+c\|q\|_{1}^{2}+c|T|_{\xi=1}|_{4}^{4}\\
&&+c|q|_{\xi=1}|_{4}^{4}+c(|Q_{1}|_{\xi=1}|_{2}^{2}+|Q_{2}|_{\xi=1}|_{2}^{2} )+c(|Q_{1}|_{2}^{2}+|Q_{2}|_{2}^{2} ),
\end{eqnarray*}
\begin{eqnarray*}
&&\frac{1}{2}\frac{d |\nabla T|_{2}^{2}}{dt}+|\Delta T|_{2}^{2}+|\nabla T_{\xi}|_{2}^{2}+\alpha_{s}|\nabla T|_{\xi=1} |_{2}^{2}\\
&\leq& c|\nabla T|_{2}^{2}
+c|Q_{1}|_{2}^{2}+c(1+|\nabla T_{\xi}|_{2})(1  +  |\Delta v |_{2} )\\
&&+c(1+|\nabla q|_{2})(1 +   |\Delta v |_{2} )\\
&\leq& \varepsilon |\nabla T_{\xi}|_{2}^{2}+c|\nabla T|_{2}^{2}
+c(1  +  |\Delta v |_{2}^{2} )+c(1+|\nabla q|_{2}^{2})+c|Q_{1}|_{2}^{2},
\end{eqnarray*}
and
\begin{eqnarray*}
&&\frac{1}{2}\frac{d |\nabla q|_{2}^{2}  }{dt}+|\Delta q|_{2}^{2}+|\nabla q_{\xi}|_{2}^{2}+\beta_{s}| \nabla q|_{\xi=1} |_{2}^{2}\\
&\leq& \varepsilon|\Delta q|_{2}^{2}+ \varepsilon |\nabla T_{\xi}|_{2}^{2}   +c|Q_{2}|_{2}^{2}+c|\nabla q|_{2}^{2}\\
&&+c(1+|\Delta v|_{2}^{2}  ).
\end{eqnarray*}
Since it is proved in $\cite{ZG}$ that $|T|_{\xi=1}|_{4}$ and $|q|_{\xi=1}|_{4} $ are uniformly bounded with respect to $t,$ summing the above estimates about $T$ and $q$ and integrating  with  respect to $t$ over $[0, \tau]$ yields
\begin{eqnarray}
&&\int_{0}^{\tau}(|\nabla T_{\xi}(t)|_{2}^{2}+|T_{\xi\xi}(t)|_{2}^{2}+|\Delta T(t)|_{2}^{2}+\alpha_{s}|\nabla T(t)|_{\xi=1}|_{2}^{2})dt\nonumber\\
&&+\int_{0}^{\tau}(|\nabla q_{\xi}(t)|_{2}^{2}+|q_{\xi\xi}(t)|_{2}^{2}+|\Delta q(t)|_{2}^{2}+\beta_{s}|\nabla q(t)|_{\xi=1}|_{2}^{2})dt\nonumber\\
&\leq&c\int_{0}^{\tau}(1+\|T(t)\|_{1}^{2}+\|q(t)\|_{1}^{2}+\|v(t)\|_{2}^{2}+|T|_{\xi=1}(t)|_{4}^{4}+ |q|_{\xi=1}(t)|_{4}^{4})dt\nonumber\\
&&+\tau(|Q_{1}|_{\xi=1}|_{2}^{2}+|Q_{2}|_{\xi=1}|_{2}^{2}+|Q_{1}|_{2}^{2}+|Q_{2}|_{2}^{2}  )\nonumber\\
&&+|\|T(\tau)\|_{1}^{2}-\|T(0)\|_{2}^{2}|+|\|q(\tau)\|_{1}^{2}-\|q(0)\|_{2}^{2}|\nonumber\\
&\leq&c(\tau^{\frac{1}{2}}+  \tau).
\end{eqnarray}
Combining $(4.45)-(4.47),$ we prove the last inequality of the theorem.
\hspace{\fill}$\square$
\end{proof}
Next, we state our main result of this paper.
\begin{theorem}
Assume $Q_{1}, Q_{2}\in L^{2}(\mho)$ and $ Q_{1}|_{\xi=1}, Q_{2}|_{\xi=1}\in L^{2}(S^{2}).$ Then the global attractor $\mathcal{A}$ has finite Hausdorff and fractal dimensions measured in the $V$ space.
\end{theorem}
\begin{proof}
Recalling that in $\cite{ZG}$ we prove that the strong solution $U$ to $(1.11)-(1.17)$ is Lipschitz continuous with respect to initial data in space $V.$ Therefore, in order
to prove our theorem, according to Theorem $2.2$, we only need to show the condition $(ii)$ of Theorem $2.2.$
\par
Let initial data $U_{i}=(v_{i}, T_{i},  q_{i} ), i=1,2$ be two strong solutions to $(1.11)-(1.17)$ with initial data  $U_{i}(0)=(v_{i,0}, T_{i,0},  q_{i,0} )  \in \mathcal{A},$ where $\mathcal{A}$ is the global attractor given by Theorem $2.1.$ In view of the invariance property of global attractor and Proposition 2.1, we can assume $ U_{i}(0) \in D(A_{1})\times D(A_{2})\times D(A_{3})$ with $i=1,2.$ By virtue of energy estimates we can easily show that $|\partial_{t}U_{i}(0)|_{2} $ are finite for $i=1,2.$ Consequently, according to Theorem $3.1,$ we know that $|\partial_{t}U_{i}(t)|$ are uniformly bounded with respect to $t(t\geq0)$. Meanwhile, we also have $U(t)\in D(A_{1})\times D(A_{2})\times D(A_{3})$ for any $t\geq 0$ in view of Theorem $2.1.$ Let's define
\begin{eqnarray*}
\mu=v^{1}-v^{2},\ \ \vartheta=T^{1}-T^{2},\ \ \rho= q^{1}-q^{2},\ \ \Phi(\theta, \varphi,t)=\Phi^{1}_{s}(\theta, \varphi,t)-\Phi^{2}_{s}(\theta, \varphi,t).
\end{eqnarray*}
Then we derive from $(1.11)-(1.17)$ that $\mu, \vartheta$ and $\rho$ satisfy
\begin{eqnarray}
\partial_{t}\mu&+&L_{1}\mu+\nabla_{v_{1}}\mu+\nabla_{\mu}v_{2}+\Big{(}\int_{\xi}^{1}\mathrm{div} v_{1}(x,y,\xi',t )d\xi'   \Big{)}\partial_{\xi}\mu\nonumber\\
&&+\Big{(}\int_{\xi}^{1}\mathrm{div} \mu(x,y,\xi',t )d\xi'   \Big{)}\partial_{\xi}v_{2}+\frac{f}{R_{0}}\mu^{\bot}+ \mathrm{grad} \Phi_{s}\nonumber\\
&&+\int_{\xi}^{1}\frac{bP}{p} \mathrm{grad} \vartheta d\xi'+\int_{\xi}^{1}\frac{abP}{p} \mathrm{grad} (q_{1} \vartheta) d\xi'+\int_{\xi}^{1}\frac{abP}{p} \mathrm{grad} (\rho T_{2}) d\xi'=0,\\
\partial_{t}\vartheta&+&L_{2}\vartheta+\nabla_{v_{1}}\vartheta+\nabla_{\mu}T_{2}+\Big{(}\int_{\xi}^{1}\mathrm{div} v_{1}(x,y,\xi',t )d\xi'   \Big{)}\partial_{\xi}\vartheta\nonumber\\
&&+\Big{(}\int_{\xi}^{1}\mathrm{div} \mu(x,y,\xi',t )d\xi'   \Big{)}\partial_{\xi}T_{2}-\frac{bP}{p}\Big{(}\int_{\xi}^{1}\mathrm{div} \mu(x,y,\xi',t )d\xi'   \Big{)}\nonumber\\
&&-\frac{abP}{p}q_{1}\Big{(}\int_{\xi}^{1}\mathrm{div} \mu(x,y,\xi',t )d\xi'   \Big{)}
-\frac{abP}{p}\rho\Big{(}\int_{\xi}^{1}\mathrm{div} v_{2}(x,y,\xi',t )d\xi'   \Big{)}=0,\\
\partial_{t}\rho&+&L_{3}\rho+\nabla_{v_{1}}\rho+\nabla_{\mu}q_{2}+\Big{(}\int_{\xi}^{1}\mathrm{div} v_{1}(x,y,\xi',t )d\xi'   \Big{)}\partial_{\xi}\rho \nonumber\\ &&+\Big{(}\int_{\xi}^{1}\mathrm{div} \mu(x,y,\xi',t )d\xi'   \Big{)}\partial_{\xi}q_{2}=0.
\end{eqnarray}
Since operators $A_{i}, i=1,2,3$ are positive selfadjoint with compact resolvent, we denote by $\{\lambda_{i,k}\}_{k\geq 1}$ the corresponding eigenvalues of $A_{i}.$ Obviously $0\leq\lambda_{i,1}<\lambda_{i,2}< \lambda_{i,3} \cdots$ and $\lambda_{i,k}\rightarrow \infty$ as $k\rightarrow\infty$ for arbitrary $i\in \{1,2,3\}.$
Let $P_{i,n},i=1,2,3$ be the orthogonal projector in $H_{i}$ onto the subspace spanned by first $n$ eigenvectors of $A_{i}.$ Denote $Q_{i,n}=I-P_{i,n},i=1,2,3$ and $n\in \mathbb{N}.$
\par
Taking inner product of $(4.48)$ with $A_{1}Q_{1,n}\mu$ in $L^{2}(\mho)$, we reach
\begin{eqnarray}
&&\frac{1}{2}\frac{d}{dt}\|Q_{1,n}\mu\|_{1}^{2}+|A_{1}Q_{1,n}\mu|_{2}^{2}\nonumber\\
&=& -\langle \nabla_{v_{1}}\mu,     A_{1}Q_{1,n}\mu \rangle-\langle \nabla_{\mu}v_{2},    A_{1}Q_{1,n}\mu \rangle\nonumber\\
&&-\langle\Big{(}\int_{\xi}^{1}\mathrm{div} \mu(x,y,\xi',t )d\xi'   \Big{)}\partial_{\xi}v_{2},    A_{1} Q_{1,n}\mu  \rangle\nonumber\\
&&-\langle \Big{(}\int_{\xi}^{1}\mathrm{div} v_{1}(x,y,\xi',t )d\xi'   \Big{)}\partial_{\xi} \mu,   A_{1}Q_{1,n}\mu \rangle\nonumber\\
&&-\langle \int_{\xi}^{1}\frac{bP}{p} \mathrm{grad }\vartheta d\xi',    A_{1} Q_{1,n}\mu \rangle-\langle (\frac{f}{R_{0}} \mu^{\bot}+\mathrm{grad} \Phi_{s}),   A_{1}Q_{1,n}\mu   \rangle\nonumber\\
&&-\langle \int_{\xi}^{1}\frac{abP}{p} \mathrm{grad}(q_{1}\vartheta)d\xi',     A_{1} Q_{1,n}\mu   \rangle-\langle \int_{\xi}^{1}\frac{abP}{p} \mathrm{grad}(\rho T_{2})d\xi'   ,     A_{1} Q_{1,n}\mu\rangle\nonumber\\
&=:&\sum_{m=1}^{8}I_{m}.
\end{eqnarray}
By  H$\mathrm{\ddot{o}}$lder inequality, Agmon inequality and Young's inequality , we have
\begin{eqnarray*}
I_{1}&\leq& |v_{1}|_{\infty}(|\nabla_{e_{\theta}}\mu|_{2}+  |\nabla_{e_{\varphi}}\mu|_{2})|A_{1}Q_{1,n}\mu|_{2}\\
&\leq &c\|v_{1}\|_{1}^{\frac{1}{2}}|A_{1}v_{1}|_{2}^{\frac{1}{2}}\|\mu\|_{1}|A_{1}Q_{1,n}\mu|_{2}\\
&\leq & \varepsilon |A_{1}Q_{1,n}\mu|_{2}^{2}+c\|v_{1}\|_{1}|A_{1}v_{1}|_{2}\|\mu\|_{1}^{2}.
\end{eqnarray*}
Similarly, we obtain
\begin{eqnarray*}
I_{2}\leq |\mu|_{\infty}\|v_{2}\|_{1}|A_{1}Q_{1,n}\mu|_{2}\leq \varepsilon |A_{1}Q_{1,n}\mu|_{2}^{2}+c\|\mu\|_{1} \|\mu\|_{2}  \|v_{2}\|_{1}^{2}.
\end{eqnarray*}
In view of Lemma $2.4$ and Young's inequality, we reach
\begin{eqnarray*}
I_{3}+I_{4}&\leq& c|A_{1}Q_{1,n}\mu|_{2}\|\mu\|_{1}^{\frac{1}{2}}\|\mu\|_{2}^{\frac{1}{2}}\|v_{2}\|_{1}^{\frac{1}{2}}\|v_{2}\|_{2}^{\frac{1}{2}}\\
&&+c|A_{1}Q_{1,n}\mu|_{2}\|\mu\|_{1}^{\frac{1}{2}}\|\mu\|_{2}^{\frac{1}{2}}\|v_{1}\|_{1}^{\frac{1}{2}}\|v_{1}\|_{2}^{\frac{1}{2}}\\
&\leq&\varepsilon |A_{1}Q_{1,n}\mu|_{2}^{2}+c\|\mu\|_{1}\|\mu\|_{2}\|v_{1}\|_{1}\|v_{1}\|_{2}\\
&&+c\|\mu\|_{1}\|\mu\|_{2}\|v_{2}\|_{1}\|v_{2}\|_{2}.
\end{eqnarray*}
Similarly, we have the estimates of $I_{7}$ and $I_{8}$ that
\begin{eqnarray*}
I_{7}+I_{8}&\leq&\varepsilon |A_{1}Q_{1,n}\mu|_{2}^{2}+ c|q_{1}|_{2}\|q_{1}\|_{1}\|\vartheta\|_{1}\|\vartheta\|_{2}\\
&&+c\|q_{1}\|_{1}\|q_{1}\|_{2}|\vartheta|_{2}\|\vartheta\|_{1}\\
&&+c|T_{2}|_{2}\|T_{2}\|_{1}\|\rho\|_{1}\|\rho\|_{2}\\
&&+c\|T_{2}\|_{1}\|T_{2}\|_{2}|\rho|_{2}\|\rho\|_{1}\\
&\leq&\varepsilon |A_{1}Q_{1,n}\mu|_{2}^{2}+c\|q_{1}\|_{1}\|q_{1}\|_{2}\|\vartheta\|_{1}\|\vartheta\|_{2}\\
&&+c\|T_{2}\|_{1}\|T_{2}\|_{2}\|\rho\|_{1}\|\rho\|_{2},
\end{eqnarray*}
where the last inequality follows by Sobolev inequality. Estimates of $I_{5}$ and $I_{6}$ follows by H$\mathrm{\ddot{o}}$lder inequality and Young's inequality
\begin{eqnarray*}
I_{5}+I_{6}\leq \varepsilon |A_{1}Q_{1,n}\mu|_{2}^{2}+c\|\vartheta\|_{1}^{2}+c|\mu|_{2}^{2}.
\end{eqnarray*}
Combining $(4.51)$ and the estimates of $I_{1}-I_{8},$ we have
\begin{eqnarray}
\frac{d}{dt}\|Q_{1,n}\mu\|_{1}^{2}+|A_{1}Q_{1,n}\mu|_{2}^{2}&\leq& c\|v_{1}\|_{1}\|v_{1}\|_{2}\|\mu\|_{1}^{2}+c\|\mu\|_{1} \|\mu\|_{2}  \|v_{2}\|_{1}^{2}\nonumber\\
&&+c\|\mu\|_{1}\|\mu\|_{2}\|v_{1}\|_{1}\|v_{1}\|_{2}\nonumber\\
&&+c\|\mu\|_{1}\|\mu\|_{2}\|v_{2}\|_{1}\|v_{2}\|_{2}\nonumber\\
&&+c\|q_{1}\|_{1}\|q_{1}\|_{2}\|\vartheta\|_{1}\|\vartheta\|_{2}\nonumber\\
&&+c\|T_{2}\|_{1}\|T_{2}\|_{2}\|\rho\|_{1}\|\rho\|_{2}\nonumber\\
&&+c\|\vartheta\|_{1}^{2}+c\|\mu\|_{1}^{2}.
\end{eqnarray}
Taking inner product of $(4.49)$ with $A_{2}Q_{2,n}\vartheta$ in $L^{2}(\mho)$, we reach
\begin{eqnarray}
&&\frac{1}{2}\frac{d}{dt}\|Q_{2,n}\vartheta\|_{1}^{2}+|A_{2}Q_{2,n}\vartheta|_{2}^{2}\nonumber\\
&=& -\langle   \nabla_{v_{1}}\vartheta,    A_{2}Q_{2,n}\vartheta\rangle  - \langle \nabla_{\mu}T_{2},    A_{2}Q_{2,n}\vartheta\rangle\nonumber\\
&&-\langle \Big{(}\int_{\xi}^{1}\mathrm{div} v_{1}(x,y,\xi',t )d\xi'   \Big{)}\partial_{\xi}\vartheta,    A_{2}Q_{2,n}\vartheta\rangle\nonumber\\
&&-\langle \Big{(}\int_{\xi}^{1}\mathrm{div} \mu(x,y,\xi',t )d\xi'   \Big{)}\partial_{\xi}T_{2},  A_{2}Q_{2,n}\vartheta\rangle\nonumber\\
&&+\langle \frac{bP}{p}\Big{(}\int_{\xi}^{1}\mathrm{div} \mu(x,y,\xi',t )d\xi'   \Big{)}, A_{2}Q_{2,n}\vartheta\rangle\nonumber\\
&&+\langle \frac{abP}{p}q_{1}\Big{(}\int_{\xi}^{1}\mathrm{div} \mu(x,y,\xi',t )d\xi'   \Big{)}, A_{2}Q_{2,n}\vartheta\rangle\nonumber\\
&&+\langle \frac{abP}{p}\rho\Big{(}\int_{\xi}^{1}\mathrm{div} v_{2}(x,y,\xi',t )d\xi'   \Big{)}, A_{2}Q_{2,n}\vartheta\rangle.
\end{eqnarray}
Obviously, $(4.53)$ is very similar to $(4.51),$ so taking an analogous argument as $(4.52)$ we have
\begin{eqnarray}
\frac{d}{dt}\|Q_{2,n}\vartheta\|_{1}^{2}+|A_{2}Q_{2,n}\vartheta|_{2}^{2}&\leq& c\|v_{1}\|_{1}\|v_{1}\|_{2}\|\vartheta\|_{1}^{2}
+c\|\mu\|_{1}\|\mu\|_{2}\|T_{2}\|_{1}^{2}\nonumber\\
&&+c\|v_{1}\|_{1}\|v_{1}\|_{2}\|\vartheta\|_{1}\|\vartheta\|_{2}\nonumber\\
&&+c\|\mu\|_{1}\|\mu\|_{2}\|T_{2}\|_{1}\|T_{2}\|_{2}\nonumber\\
&&+c\|\mu\|_{1}\|\mu\|_{2}|q_{1}|_{2}\|q_{1}\|_{1}\nonumber\\
&&+c\|v_{2}\|_{1}\|v_{2}\|_{2}|\rho|_{2}\|\rho\|_{1}+c\|\mu\|_{1}^{2}.
\end{eqnarray}
Similarly,  we obtain the estimates of $Q_{3,n}\rho:$
\begin{eqnarray}
\frac{d}{dt}\|Q_{3,n}\rho\|_{1}^{2}+|A_{3}Q_{3,n}\rho|_{2}^{2}&\leq& c\|v_{1}\|_{1}\|v_{1}\|_{2}\|\rho\|_{1}^{2}+ c\|\mu\|_{1}\|\mu\|_{2}\|q_{2}\|_{1}^{2}\nonumber\\
&&+c\|v_{1}\|_{1}\|v_{1}\|_{2}\|\rho\|_{1}\|\rho\|_{2}+c\|\mu\|_{1}\|\mu\|_{2}\|q_{2}\|_{1}\|q_{2}\|_{2}.
\end{eqnarray}
For $t\geq 0,$ let
\begin{eqnarray*}
&&\varphi(t):=\|Q_{1,n}\mu(t)\|_{1}^{2}+\|Q_{2,n}\vartheta(t)\|_{1}^{2}+\|Q_{3,n}\rho(t)\|_{1}^{2},\\
&&\psi(t):=\|\mu(t)\|_{1}^{2}+\|\vartheta(t)\|_{1}^{2}+\|\rho(t)\|_{1}^{2}.
\end{eqnarray*}
Summing $(4.52), (4.54) $ and $(4.55)$ and using the uniform-time estimates of the strong solution to $(1.11)-(1.17)$ in space $V$ we get
\begin{eqnarray}
\frac{d}{dt}\varphi+\lambda_{n}\varphi&\leq&c\psi(1+\|v_{1}\|_{2}+\|v_{2}\|_{2} )\nonumber\\
&&+c\|\mu\|_{1}\|\mu\|_{2}(1+ \|v_{1}\|_{2}+\|v_{2}\|_{2}+\|T_{2}\|_{2}+\|q_{2}\|_{2}   )\nonumber\\
&&+c\|\vartheta\|_{1}\|\vartheta\|_{2}(\|q_{1}\|_{2}+\|v_{1}\|_{2})\nonumber\\
&&+c\|\rho\|_{1}\|\rho\|_{2}(\|T_{2}\|_{2}+\|v_{1}\|_{2}+\|v_{2}\|_{2}  ),
\end{eqnarray}
where $\lambda_{n}:=min\{ \lambda_{1,n},  \lambda_{2,n}, \lambda_{3,n}\}.$ Obviously $0\leq \lambda_{n}\rightarrow \infty $ as $n\rightarrow \infty. $
Integration $(4.56)$ with respect to $t$ over $[0,T]$ yields
\begin{eqnarray}
\varphi(T)&\leq& e^{-\lambda_{n}T}\varphi(0)+ ce^{-\lambda_{n}T}\int_{0}^{T}e^{\lambda_{n}t}\psi(t)dt\nonumber\\
&&+ce^{-\lambda_{n}T}\int_{0}^{T}e^{\lambda_{n}t}\psi(t)(\|v_{1}(t)\|_{2}+\|v_{2}(t)\|_{2} )dt\nonumber\\
&&+ce^{-\lambda_{n}T}\int_{0}^{T}e^{\lambda_{n}t}\psi(t)^{\frac{1}{2}}( \|\mu(t)\|_{2}+\|\vartheta(t)\|_{2}+ \|\rho(t)\|_{2}  )\nonumber\\
&&\cdot( \|v_{1}(t)\|_{2}+\|v_{2}(t)\|_{2}+ \|T_{2}(t)\|_{2}+\|q_{1}(t)\|_{2}+ \|q_{2}(t)\|_{2} )dt\nonumber\\
&=:&\sum_{m=1}^{4}J_{m}.
\end{eqnarray}
To verify  that $\varphi(T)$ satisfies the condition $(ii)$ of Theorem $2.2,$ we need an estimate derived from $\cite{ZG}$ that
\begin{eqnarray}
\psi(t)+\int_{0}^{t}( \|\mu(t)\|_{2}^{2}+\|\vartheta(t)\|_{2}^{2}+ \|\rho(t)\|_{2}^{2}  )dt\leq   \gamma(t)\psi(0),
\end{eqnarray}
where $\gamma(t), t\geq 0, $ is a positive continuous and non-decreasing function independent of initial data.
Therefore, the estimates of $J_{2}$ follows that
\begin{eqnarray*}
J_{2}\leq c\gamma(T)\psi(0)e^{-\lambda_{n}T}\int_{0}^{T}e^{\lambda_{n}t}dt\leq c\gamma(T)\psi(0)\lambda_{n}^{-1}.
\end{eqnarray*}
Using $(4.58)$ again and Theorem $4.1$ we have
\begin{eqnarray*}
J_{3}&\leq& c\gamma(T)\psi(0)e^{-\lambda_{n}T}\int_{0}^{T}e^{\lambda_{n}t}(\|v_{1}\|_{2}+\|v_{2}\|_{2})dt\nonumber\\
&\leq& c\gamma(T)\psi(0)\int_{0}^{T} (\|v_{1}\|_{2}+\|v_{2}\|_{2})dt\nonumber\\
&\leq& c\gamma(T)\psi(0)T^{\frac{1}{2}}(\int_{0}^{T} (\|v_{1}\|_{2}^{2}+\|v_{2}\|_{2}^{2})dt )^{\frac{1}{2}}\nonumber\\
&\leq& c\gamma(T)\psi(0)(T+T^{\frac{3}{4}} ).
\end{eqnarray*}
Similarly, by $(4.58)$,  Theorem $4.1$ and H$\mathrm{\ddot{o}}$lder inequality we obtain the estimates of $J_{4}$
\begin{eqnarray*}
J_{4}&\leq& c\gamma(T)^{\frac{1}{2}}\psi(0)^{\frac{1}{2}}\Big{(}\int_{0}^{T}( \|\mu\|_{2}^{2}+\|\vartheta\|_{2}^{2}+  \|\rho\|_{2}^{2})dt \Big{)}^{\frac{1}{2}}\nonumber\\
&&\cdot \Big{(}\int_{0}^{T}( \|q_{1}\|_{2}^{2}+\|v_{1}\|_{2}^{2}+  \|v_{2}\|_{2}^{2}+  \|T_{2}\|_{2}^{2})dt \Big{)}^{\frac{1}{2}}\nonumber\\
&\leq& c\gamma(T)^{\frac{1}{2}}\psi(0)^{\frac{1}{2}}\gamma(T)^{\frac{1}{2}}\psi(0)^{\frac{1}{2}}(T^{\frac{1}{2}} +T)^{\frac{1}{2}}\nonumber\\
&\leq& c\gamma(T)\psi(0)(T^{\frac{1}{2}} +T)^{\frac{1}{2}}.
\end{eqnarray*}
Summing $(4.57)$ and estimates of $J_{1}-J_{4}$, we get
\begin{eqnarray*}
\varphi(T)\leq c\psi(0)[e^{-\lambda_{n}T}+  \gamma(T)\lambda_{n}^{-1}+ \gamma(T)(T+T^{\frac{3}{4}})+  \gamma(T)(T^{\frac{1}{4}}+T^{\frac{1}{2}}) ].
\end{eqnarray*}
Since $U_{i}(0)\in \mathcal{A} $ and $\mathcal{A}$ is bounded in $V$, $ \psi(0)$ is bounded. We can choose $T$ which is small enough such that for $\delta\in (0,1)$
\begin{eqnarray*}
c\gamma(T)(T+T^{\frac{3}{4}})+  c\gamma(T)(T^{\frac{1}{4}}+T^{\frac{1}{2}})\leq \frac{\delta}{2}.
\end{eqnarray*}
Then let $\lambda_{n}$ be big enough such that
\begin{eqnarray*}
ce^{-\lambda_{n}T}+c\gamma(T)\lambda_{n}^{-1}\leq \frac{\delta}{2}.
\end{eqnarray*}
Consequently,  for small $T$ and large $\lambda_{n},$ we have
\begin{eqnarray*}
\varphi(T)\leq \delta \psi(0),
\end{eqnarray*}
which combined with Theorem 2.2 proves this theorem.
\hspace{\fill}$\square$
\end{proof}

\end{theorem}

\def\refname{ Bibliography}

\end{document}